\documentclass[12pt]{article}

\usepackage{tikz}
\usetikzlibrary{arrows.meta,calc,positioning,decorations.pathreplacing}

\usepackage{graphicx}
\usepackage{caption}
\usepackage{amsfonts}
\usepackage{amscd}
\usepackage{indentfirst}
\usepackage{fancyhdr,epsfig,color,mathtools}
\usepackage{amsmath,amssymb,amsthm, amsgen, amstext, amsbsy, amsopn,yfonts,bbm}
\usepackage{xcolor}
\usepackage{mathrsfs}

\usepackage{tocloft}

\cftsetindents{section}{5em}{3em}
\cftsetindents{subsection}{5em}{3em}

\usepackage[bookmarks=true, bookmarksnumbered=true, bookmarksopen=true, backref=page]{hyperref}
\renewcommand*{\backrefalt}[4]{%
	\ifcase #1 (Not cited.)%
	\or        (Cited on page~#2.)%
	\else      (Cited on pages~#2.)%
	\fi}
\hypersetup{
	colorlinks   = true,
	citecolor    = magenta,
	linkcolor    =blue,
	urlcolor     =magenta	
}

\newtheorem{thm}{Theorem}[section]

\newtheorem{prop}[thm]{Proposition}
\newtheorem{cor}[thm]{Corollary}
\newtheorem{defin}[thm]{Definition}
\newtheorem{rem}[thm]{Remark}
\newtheorem{exam}[thm]{Example}

\newtheoremstyle{condstyle} 
  {3pt}                         
  {3pt}                         
  {\slshape}                    
  {}                            
  {\scshape}                    
  {.}                           
  {.5em}                        
  {}                            

\theoremstyle{condstyle}

\newcommand{\F}{{\mathbb{F}}}

\newcommand{\I}{{\mathbb{I}}}
\newcommand{\Id}{{\textit{Id}}}
\newcommand{\J}{{\mathbb{J}}}

\newcommand{\R}{{\mathbb{R}}}
\newcommand{\SP}{{\mathbb{S}}}

\newcommand{\Z}{{\mathbb{Z}}}

\newcommand{\cA}{{\mathcal{A}}}
\newcommand{\cB}{{\mathcal{B}}}

\newcommand{\cI}{{\mathcal{I}}}
\newcommand{\cJ}{{\mathcal{J}}}

\newcommand{\cM}{{\mathcal{M}}}

\newcommand{\cR}{{\mathcal{R}}}

\newcommand{\bL}{{\bar{L}}}

\newcommand{\bY}{{\bar{Y}}}

\newcommand{\barv}{{\bar{v}}}

\newcommand{\btheta}{\bar{\theta}}
\newcommand{\bvartheta}{\bar{\vartheta}}

\newcommand{\tc}{{\widetilde{c}}}

\newcommand{\tI}{{\widetilde{I}}}

\newcommand{\tL}{{\widetilde{L}}}
\newcommand{\tM}{{\widetilde{M}}}

\newcommand{\tV}{{\widetilde{V}}}

\newcommand{\tW}{{\widetilde{W}}}

\newcommand{\tpi}{\widetilde{\pi}}

\newcommand{\tlambda}{\widetilde{\lambda}}

\newcommand{\tLambda}{\widetilde{\Lambda}}

\newcommand{\tPhi}{\widetilde{\Phi}}

\newcommand{\ttheta}{\widetilde{\theta}}

\newcommand{\lmininfty}{{\it l_{min, [\infty]}}}

\newcommand{\lminc}{{\it l_{min,[c]}}}
\newcommand{\lmintc}{{\it l_{min,[\tc]}}}

\newcommand{\Ups}{\Upsilon}

\newcommand{\SSigma}{{{\mathcal S}\Sigma}}

\newcommand{\df}{\textit{def}\,}

\newcommand{\graph}{{\it graph\,}}

\newcommand{\Image}{{\rm Im\,}}

\newcommand{\Ker}{{\rm Ker\,}}

\newcommand{\Qed}{\ \hfill \qedsymbol \bigskip}

\newcommand{\Sh}{{\it Sh}}
\newcommand{\spec}{{\it spec}}

\newcommand{\Op}{{\it Op}}

\newcommand{\pr}{\hbox{\textrm pr}}

\newcommand{\sgrad}{{\hbox{\rm sgrad\,}}}

\newcommand{\gap}{G}

\begin{document}

\title{Lagrangian cobordisms and morphisms of persistence modules}

\author{%
\textsc{Michael Entov}$^{1}$,\ \textsc{Leonid Polterovich}$^{2}$\\
\ \\
\textit{To the memory of Vladimir I. Arnol'd, with gratitude}
\ \\
\
}

\date{\today}

\maketitle

\footnotetext[1]{Partially supported by the Israel Science
Foundation grant \# 2033/23.
}

\footnotetext[2]{Partially supported by the Israel Science Foundation - NSFC grant \#3231/23.
}

\begin{abstract}
We consider a version of Legendrian contact homology filtered by the action of Reeb chords.
In this version, a union of two disjoint exact Lagrangian cobordisms with Legendrian boundary conditions gives rise to a morphism between
persistence modules. We discuss a certain algebraic invariant of such morphisms and use it to obtain
obstructions to Lagrangian isotopies of the cobordisms, as well as applications to symplectic and contact dynamics.
\end{abstract}

\tableofcontents

\section{Introduction and a sample application}
\label{sec-intro-and-sample-application}

In this paper we consider Lagrangian cobordisms and their isotopies, basic
objects of symplectic and contact topology. One of the principal founders of this field was V. I.
Arnol’d, to whose memory this paper is dedicated.


By a Lagrangian cobordism we mean
a compact Lagrangian submanifold $L$ with positive and negative Legendrian boundaries; $L$ is supposed to lie in
an
ambient symplectic cobordism $(M,\theta)$ between contact manifolds and the positive/negative boundary of $L$
lies, respectively, in the positive/negative boundary of $M$.
(Here and further on,
$\theta$ is a 1-form such
that $d\theta$ is a symplectic form).
We also assume that the Lagrangian cobordism is cylindrical near its boundaries.
Such a cobordism $L\subset (M,\theta)$ is called {\it exact (with respect to $\theta$)} if $\theta|_L = df$ for
a smooth function $f:L\to\R$ that vanishes near the negative boundary of $L$ and is
identically equal to a constant $A_{L,\theta}$ near the positive boundary of $L$.

In fact, we will consider Lagrangian cobordisms of the form $L=L_0\sqcup L_1$ which are
unions of ordered pairs of disjoint
(possibly disconnected) exact Lagrangian cobordisms; we will call these
{\it two-part-exact Lagrangian cobordisms}.
 The positive and negative
boundaries of such Lagrangian cobordisms
are the so-called {\it two-part Legendrian submanifolds} -- i.e. the
unions of ordered pairs of compact (possibly disconnected)
disjoint Legendrian submanifolds.
The {\it gap} of a two-part-exact Lagrangian cobordism $L=L_0\sqcup L_1\subset (M,\theta)$ (with respect to $\theta$)
is defined as
\[
\gap_{L,\theta} := A_{L_1,\theta} - A_{L_0,\theta}.
\]

The main object discussed in this paper is a new quantity associated to two-part-exact Lagrangian cobordisms.
Roughly speaking, to a pair $(\Lambda,\lambda)$, formed by a
two-part Legendrian submanifold $\Lambda_0\sqcup\Lambda_1$ in $\Sigma$
and a contact form $\lambda$ on $\Sigma$,
we associate
a persistence module obtained from a version of the Legendrian contact homology
of $\Lambda_0\sqcup\Lambda_1$, with respect to $\lambda$,
filtered by the action.
To a two-part-exact Lagrangian cobordism
between the pairs as above we associate a morphism between the corresponding
persistence modules obtained using the morphism between the (non-filtered) Legendrian contact homologies
constructed in the works of Ekholm \cite{Ekh-JEMS} and Ekholm-Honda-K\'alm\'an \cite{Ekholm-Honda-Kalman}.
Such morphisms of persistence modules are rich algebraic objects.
We define an algebraic invariant of these morphisms that behaves ``nicely'' under composition of morphisms,
and consequently under concatenation of the cobordisms. It is also invariant under Lagrangian
isotopies of the cobordisms that are fixed near the boundary and keep the gap fixed.
This invariant gives rise to geometric constraints
on the Lagrangian isotopy classes of the cobordisms -- see Section~\ref{subsec-appl-restriction-on-gap-fcns}.
It also has applications to contact and symplectic Hamiltonian dynamics -- see Section~\ref{subsec-appls-to-dynamics}.

Here is a sample application of the first kind. Let us emphasize that our general results yield similar applications
for a much wider class of Legendrian submanifolds.

Let $Q$ be a closed connected manifold, $\dim Q >0$. Let $\Sigma := J^1 Q$ and let $\lambda$ be the standard contact 1-form
defining the standard contact structure on $\Sigma$. Let $F: Q\to\R$ be a positive Morse function and let $\kappa>0$.
Define
\[
a:= \min_Q F >0,
\]
\[
\Lambda^{\prime}:= j^1 0,
\]
\[
\Lambda^{\prime\prime} := j^1 F,
\]
\[
\Lambda^{\prime\prime\prime} := j^1 (F+\kappa).
\]

Clearly, $\Lambda^{\prime}, \Lambda^{\prime\prime}, \Lambda^{\prime\prime\prime} \subset (\Sigma,\xi)$ are
disjoint compact Legendrian submanifolds.

Consider the trivial Lagrangian cobordism
\[
L:= (\Lambda^{\prime}\sqcup\Lambda^{\prime\prime})\times [s_-,s_+]\subset (\Sigma\times [s_-,s_+],s\lambda).
\]


\hfill


\begin{thm}
\label{thm-intro-application}
Let
\[
0< s_- < s_0 < s_1 < s_+.
\]
Assume there exists an isotopy, fixed near the boundary, in the class of two-part-exact
Lagrangian cobordisms in $(\Sigma\times [s_-,s_+],s\lambda)$, from $L$ to a two-part-exact
Lagrangian cobordism $\tL = \tL_0\sqcup\tL_1\subset(\Sigma\times [s_-,s_+],s\lambda)$
whose intersection
with $\Sigma\times [s_0-\delta,s_1+\delta]$, for some $\delta>0$,
is the trivial Lagrangian cobordism
$(\Lambda^{\prime}\sqcup\Lambda^{\prime\prime\prime})\times [s_0-\delta,s_1+\delta]\subset (\Sigma\times [s_-,s_+],s\lambda)$.

Then
\begin{equation}
\label{eqn-bound-thm-intro-application}
\gap_{\tL,s\lambda} < (s_+ - s_-) a -
(s_1 - s_0) (a+\kappa).
\end{equation}


\end{thm}


\hfill


\begin{figure}[htbp]
\centering
\includegraphics[width=0.7\textwidth]{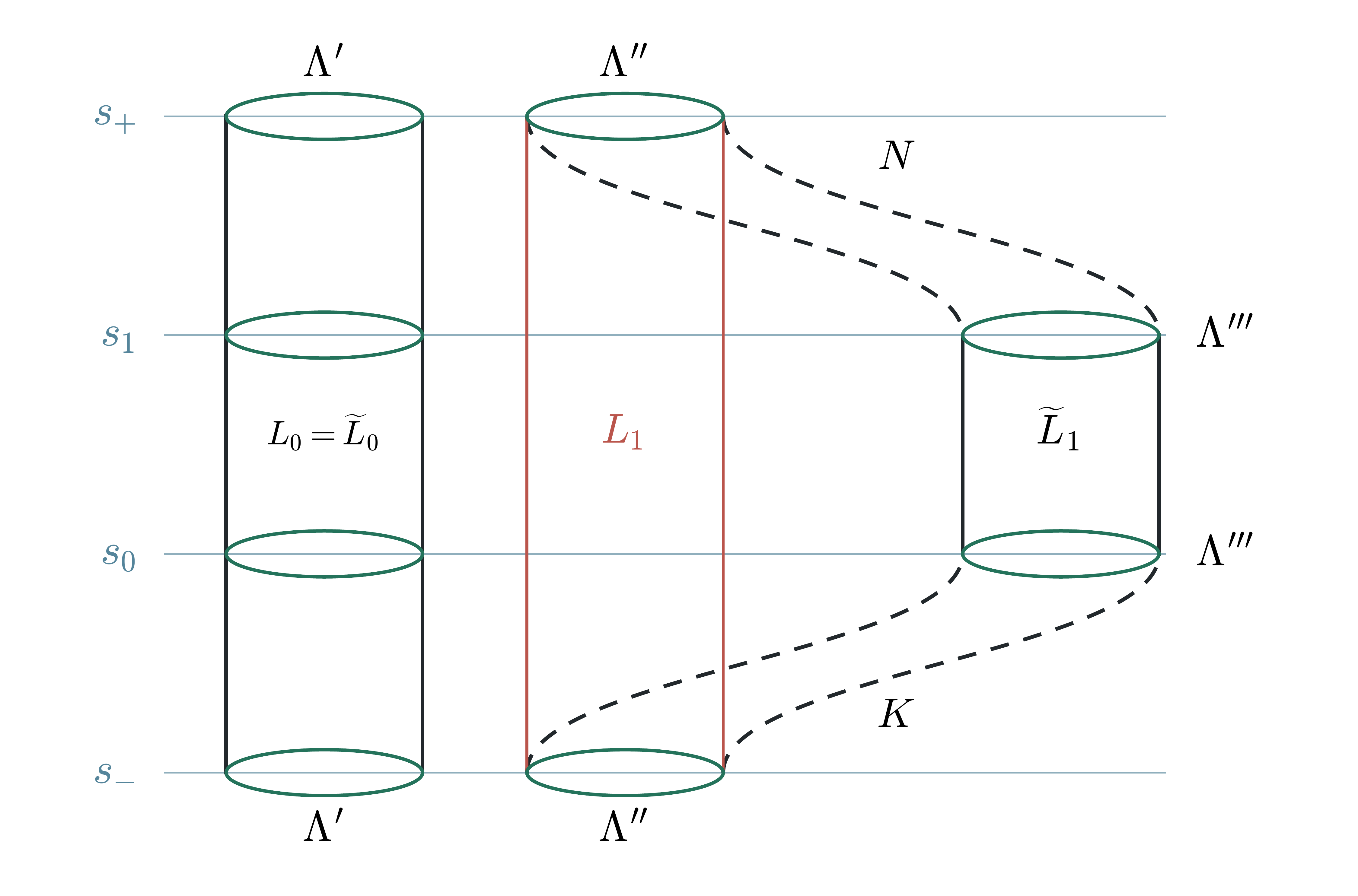}
\caption{}
\label{fig-cobordism-schematic-intro}
\end{figure}


\hfill


See Figure~\ref{fig-cobordism-schematic-intro} for a geometric illustration of the Lagrangian
cobordisms in Theorem~\ref{thm-intro-application}.

For the proof of Theorem~\ref{thm-intro-application}, see Theorem~\ref{thm-intro-application-copy}.

The bound \eqref{eqn-bound-thm-intro-application} is sharp -- for details see Example~\ref{exam-estimate-sharp}.

The non-trivial phenomenon appearing in Theorem~\ref{thm-intro-application} for $\dim Q>0$
generalizes an elementary phenomenon in the case where $Q$ and the Legendrians are $0$-dimensional,
as discussed in the example below.

Further on, $\R_+ := (0,+\infty)$.


\hfill


\begin{exam}
\label{exam-Sigma-is-R-intro}
{\rm
Assume $\Sigma=\R (z)$, $\lambda= dz$. (Note that the kernel distribution of $dz$, namely
the zero distribution of the tangent spaces of $\R$, is a contact structure $\xi$ on $\R$
and $dz$ is a contact form for $\xi$, whose Reeb vector field is $\partial/\partial z$). Let
\[
0<s_-<s_0<s_1<s_+.
\]
Consider disjoint Legendrian submanifolds
\[
\Lambda^{\prime} := \{ z_1\},\
\Lambda^{\prime\prime} := \{ z_2\},\
\Lambda^{\prime\prime\prime} := \{ z_3\}
\]
of $\R$, where
$z_1 <  z_2 < z_3$. In this case,
\[
a = z_2 - z_1,\ a+\kappa = z_3 - z_1.
\]

Set
\[
L_i:= \{ z_{i+1}\}\times [s_-,s_+]\subset \R\times [s_-,s_+],\ i=0,1.
\]
Then $L:= L_0\sqcup L_1\subset (\R\times [s_-,s_+], sdz)$
is a trivial two-part-exact Lagrangian
cobordism and $\Lambda^{\prime}\sqcup \Lambda^{\prime\prime\prime}$ is a two-part Legendrian submanifold of $(\R,\xi)$.

Consider an isotopy
$\{ L^\tau = L^\tau_0\sqcup L^\tau_1\}_{0\leq\tau\leq T}$ of $L=L^0$,
fixed near the boundary, in the class of
two-part-exact Lagrangian cobordisms in $(\R\times [s_-,s_+], sdz)$.

One easily checks that
\[
\gap_{L^\tau,sdz} = (s_+ - s_-) (z_2-z_1) - A(\tau),
\]
where $A(\tau)$ is
the $dz\wedge ds$-area of the curvilinear quadrilateral in $\R (z)\times\R_+ (s)$ bounded by
$\{s=s_-\}$, $\{ s=s_+\}$, $L^\tau_0$ and $L^\tau_1$.

Assume, in addition, that, for some small $\delta>0$,
the intersection of $L^T$ with $\R\times [s_0-\delta,s_1+\delta]$
is $(\Lambda^{\prime}\sqcup \Lambda^{\prime\prime\prime})\times [s_0-\delta,s_1+\delta]$. Then
\[
A(T)> (s_1-s_0)(z_3-z_1),
\]
hence
\[
\gap_{L^T,sdz}< (s_+ - s_-) (z_2-z_1) - (s_1-s_0)(z_3-z_1) = (s_+ - s_-) a - (s_1-s_0)(a+\kappa),
\]
which is a version of \eqref{eqn-bound-thm-intro-application}.

Let us note that one can easily construct an isotopy as above so that
$\gap_{L^T,sdz}$ is arbitrarily close to $(s_+ - s_-) (z_2-z_1) - (s_1-s_0)(z_3-z_1)$ from below.

}
\end{exam}


\hfill


Similar geometric constraints on Lagrangian cobordisms between Legendrians were studied by J.Sabloff and L.Traynor in \cite{Sabloff-Traynor} and their work is a substantial inspiration for this paper. Note however that, unlike in our setting,
in \cite{Sabloff-Traynor} a Lagrangian cobordism is a manifold {\sl without boundary},
with cylindrical ends,
and the Lagrangian isotopies live in the whole symplectization $\Sigma\times\R_+$,
and not just in
$\Sigma\times [s_-,s_+]$.
See Remark~\ref{rem-sabloff-traynor} for a more detailed discussion of the relation between Sabloff-Traynor results and
the ones in this paper.

Let us comment on the proofs.

The technology used in the proof of the main results, of which Theorem~\ref{thm-intro-application} is
a corollary, comes from a large project \cite{EP-big}
that develops the so-called {\sl Filtered Relative Symplectic Field Theory (FRSFT)}; a part of this
project has already appeared in \cite{EP-JFPTA2022}. Roughly speaking,
the FRSFT uses various versions of Legendrian contact homology filtered by the action of the Reeb chords
in order to construct a functor from a geometric category to an algebraic one. In the basic setting,
the geometric category is formed by pairs $(\Lambda,\lambda)$, where $\Lambda$ is a (two-part) Legendrian submanifold of
a contact manifold $(\Sigma,\xi)$ and $\lambda$ is a contact
form on $(\Sigma,\xi)$. The morphisms are given by exact/two-part-exact Lagrangian cobordisms between
the Legendrian submanifolds living in exact symplectic cobordisms between the corresponding contact manifolds.
The algebraic category is formed by persistence modules, whose barcodes lie in $(0,+\infty)$, and morphisms between
them.

In this paper we will outline how a particular version of the FRSFT works out if we consider a particular class
of contact manifolds $(\Sigma,\xi)$ and of contact forms $\lambda$ on $(\Sigma,\xi)$
and restrict ourselves only to two-part-exact Lagrangian cobordisms lying in exact symplectic
cobordisms of a special kind that are subsets of the symplectization of such a $(\Sigma,\xi)$.


\hfill


\section{Preliminaries}
\label{sec-prelims}

Let $(\Sigma,\xi)$ be a contact manifold with a co-oriented contact structure.

{\bf Further on, $\Sigma$ is always assumed to be connected and
all the contact forms on $(\Sigma,\xi)$ are assumed to be compatible with the co-orientation.

All Legendrian submanifolds under consideration are assumed to be compact.}

If $\lambda,\tlambda$ are contact forms on $(\Sigma,\xi)$, we write $\lambda\prec \tlambda$
if $\lambda/\tlambda\leq c$ for some $0<c<1$.

A {\it two-part Legendrian submanifold} of $(\Sigma,\xi)$ is an ordered disjoint union of two
(not necessarily connected) Legendrian submanifolds of $(\Sigma,\xi)$.


\hfill


\subsection{Intergraphical exact symplectic cobordisms}
\label{subsec-prelims-intergraphical}

Consider the symplectization $\SSigma$ of $(\Sigma,\xi)$: it is the space of linear 1-forms on the tangent spaces of $\Sigma$ that vanish on the contact planes and are compatible with the co-orientation of $\xi$; the space $\SSigma$ comes with the canonical projection $pr:\SSigma\to \Sigma$, the natural multiplicative $\R^+$-action and a canonical 1-form $\theta_{can}$ whose differential is the standard symplectic form on $\SSigma$ -- see \cite[Appendix 4]{Arnold-math-methods} for details.

For a set $\Lambda\subset\Sigma$, denote by $S\Lambda := \pr^{-1} (\Lambda)\subset\SSigma$ the cylinder over $\Lambda$.

Any contact form $\mu$ on $(\Sigma,\xi)$ defines a hypersurface $\graph \mu\subset \SSigma$. The 1-form $\theta_{can}|_{\graph \mu}$ is a contact form for a contact structure on $\graph\mu$.
The projection
\[
pr: (\graph\mu, \ker \theta_{can}|_{\graph \mu})\to (\Sigma,\xi)
\]
is a contactomorphism identifying the contact form $\theta_{can}|_{\graph \mu}$ with $\mu$. The
Legendrian submanifolds of $(\graph\mu, \ker \theta_{can}|_{\graph \mu})$ are exactly the lifts of the Legendrian submanifolds of $(\Sigma,\xi)$ under the projection.

If $\mu_0$ is a chosen contact form on $(\Sigma,\xi)$, then there is a canonical diffeomorphism $\SSigma\to \Sigma\times\R_+ (s)$ identifying $\theta_{can}$ with $s\mu_0$. This diffeomorphism identifies $\graph\mu$ (for a contact form $\mu$ on $\Sigma$) with the set
\[
\left\{ \left(x,\left(\mu/\mu_0\right) (x)\right), x\in\Sigma\right\}\subset \Sigma\times\R_+.
\]

Given contact forms $\mu^-, \mu^+$ on $(\Sigma,\xi)$, $\mu^- \prec \mu^+$,
define $\SSigma_{\mu^-,\mu^+}\subset \SSigma$ as the codimension-zero submanifold with boundary of $\SSigma$ bounded by $\graph\mu^-$ and $\graph\mu^+$.


\hfill


\begin{defin}
\label{def-integraphical-exact-symplec-cob}
{\rm
An  {\it intergraphical exact symplectic cobordism $(M,\theta)$ based on $(\Sigma,\xi)$} is a pair formed by a
manifold with boundary $M$,
of the form $M=\SSigma_{\mu^-,\mu^+}$,
and a 1-form
$\theta$, so that for a sufficiently small $\delta>0$ and a collar neighborhood $U_\delta:=\SSigma_{\mu^-,(1+\delta)\mu^-}\cup\SSigma_{(1-\delta)\mu^+,\mu^+}$ of $\partial M$,

\begin{itemize}

\item{} for some $c>0$, $\theta=c\theta_{can}$ on $U_\delta$,

\item{} the 2-form
$d\theta$
is symplectic and is
obtained from the symplectic form
$d(c\theta_{can})$
by a diffeomorphism $M\to M$
which is the identity on $U_\delta$.

\end{itemize}

We will also say that such an $(M,\theta)$ is {\it an intergraphical exact symplectic cobordism
between $\left(\Sigma^\pm := \graph \mu^\pm, \lambda^\pm := c\mu^\pm = \theta|_{\Sigma^\pm}\right)$}.

The boundary $\Sigma^+$, respectively $\Sigma^-$,
is called {\it the positive,  respectively the negative, boundary} of $M$.

A 1-form $\theta$ on $M$ as above will be called a {\it cobordism 1-form}. It restricts to
the contact forms $\lambda^\pm = c\mu^\pm$ on the positive/negative boundary of $M$.

An intergraphical exact symplectic cobordism (based on $(\Sigma,\xi)$) is called {\it trivial} if it has the form
$(M=\SSigma_{s_-\lambda,s_+\lambda},\theta = \theta_{can})$
for some $0<s_-<s_+$ and a contact form $\lambda$ on $(\Sigma,\xi)$.


}
\end{defin}


\hfill


Equivalently, under the identification above,
trivial exact symplectic co\-bor\-disms can be described as pairs $(\Sigma\times [s_-,s_+], s\lambda)$,
where $\lambda$ is a contact form on $(\Sigma,\xi)$ and $s$ is the coordinate on $[s_-,s_+]$.


\hfill


\begin{rem}
\label{rem-rescaling-integraphical-exact-sympl-cobs}
{\rm
The multiplicative action of $c\in\R_+$ on $S\Sigma$ defines a rescaling diffeomorphism
\[
(\SSigma_{\mu^-,\mu^+}, \theta_{can})\to (\SSigma_{c\mu^-,c\mu^+}, \theta_{can}/c)
\]
between the intergraphical exact symplectic cobordisms.
}
\end{rem}


\hfill


Let $(M,\theta)$ be an intergraphical exact symplectic cobordism based on $(\Sigma,\xi)$
whose positive boundary is $(\Sigma,\lambda^+)$.
Assume $(\tM,\ttheta)$ is another intergraphical exact symplectic cobordism based on $(\Sigma,\xi)$
whose negative boundary is also
$(\Sigma, \lambda^+)$.
One can glue $(M,\theta)$ and $(\tM,\ttheta)$ and the forms $\theta$, $\ttheta$
along their common boundary into a smooth manifold $M\sharp \tM$ carrying a smooth 1-form
$\theta\sharp \ttheta$ that restricts to $\theta$ on $M$ and to $\ttheta$ on $\tM$.
The resulting
$(M\sharp \tM, \theta\sharp \ttheta)$ is again an intergraphical exact symplectic cobordism
based on $(\Sigma,\xi)$,
which we will call {\it the
concatenation of $(M,\theta)$ and $(\tM,\ttheta)$}.


\hfill


\subsection{Lagrangian cobordisms}
\label{subsec-leg-submfds-lagr-cobs}

Assume that $(M,\theta)$ is an intergraphical exact symplectic cobordism, based on $(\Sigma,\xi)$,
between $(\Sigma^\pm,\lambda^\pm)$.


Let $\Lambda^\pm = \Lambda^\pm_0\sqcup\Lambda^\pm_1\subset (\Sigma,\xi)$
be two-part Legendrian submanifolds, viewed, under the identifications from Section~\ref{subsec-prelims-intergraphical},
respectively, as submanifolds of
$(\Sigma^\pm, \ker \theta|_{\Sigma^\pm})$.

A {\it Lagrangian cobordism $L$ in $(M,\theta)$ between
$(\Lambda^\pm,\lambda^\pm)$} is a smooth compact cobordism in $M$ between
$\Lambda^+\subset \Sigma^+$ and
$\Lambda^-\subset \Sigma^-$ which is a Lagrangian submanifold of $(M, d\theta)$
and is cylindrical near the boundary -- i.e. in a neighborhood of $\Sigma^\pm$, the submanifold $L$
coincides with $S\Lambda^\pm$.

The Legendrian submanifolds $\Lambda^\pm$ are called {\it the positive/negative boundaries of $L$}.

In this paper, we consider only Lagrangian cobordisms all of whose connected components
have non-empty positive and negative boundaries.

A Lagrangian cobordism $L$ as above is called {\it exact}, or {\it one-part-exact},
if $\theta|_L = df$ for a smooth function $f: L\to\R$
which is zero near $\Lambda^-$ and identically equal to a constant $A_{L,\theta}$ near $\Lambda^+$.
The constant $A_{L,\theta}$ is independent of the choice of $f$ (as long as such an $f$ exists).

A {\it two-part-exact Lagrangian cobordism in $(M,\theta)$
between two-part Legendrian submanifolds
$\Lambda^\pm=\Lambda^\pm_0\sqcup \Lambda^\pm_1$} is a disjoint union $L= L_0\sqcup L_1$, where
each $L_i\subset (M,\theta)$, $i=0,1$ is an
exact Lagrangian
cobordism between $(\Lambda^\pm_i,\lambda^\pm)$.
The cobordisms
$L_i\subset (M,\theta)$, $i=0,1$, are called {\it the parts of $L$},
or {\it the $0$- and $1$-parts of $L$}.

Recall from Section~\ref{sec-intro-and-sample-application} that the constant
\[
\gap_{L,\theta}:=A_{L_1,\theta}-A_{L_0,\theta}
\]
is called {\it the gap of $(L,\theta)$}, or {\it the gap of $L$ with respect to $\theta$}.

Note that the definition of a two-part-exact Lagrangian cobordism includes the numbering of its parts. Switching
the numbering of the two parts changes the gap to its opposite.

A {\it trivial two-part-exact Lagrangian cobordism} is a two-part-exact Lagrangian cobordism of the form
$L:=\Lambda\times [s_-,s_+]\subset (\Sigma\times [s_-,s_+],s\lambda)$ in
a trivial exact symplectic cobordism $(\Sigma\times [s_-,s_+],s\lambda)$, for some $0<s_-<s_+$ and $\lambda$,
where $\Lambda=\Lambda_0\sqcup\Lambda_1$ is a two-part Legendrian submanifold of $(\Sigma,\xi)$. In this case,
$\gap_{L,s\lambda} =0$.

Assuming that $(M,\theta)$ and $(\tM,\ttheta)$ are intergraphical exact symplectic cobordisms (based on $(\Sigma,\xi)$) whose concatenation
$(M\sharp \tM, \theta\sharp \ttheta)$ is well-defined and $L\subset (M,\theta)$, $\tL\subset (\tM,\ttheta)$
are two-part-exact Lagrangian cobordisms with a common boundary, one can glue $L$ and $\tL$ along the common boundary and
get a two-part-exact Lagrangian cobordism
of $(M\sharp \tM, \theta\sharp \ttheta)$
-- it will be denoted by $L\sharp\tL\subset (M\sharp \tM, \theta\sharp \ttheta)$.

Assume that
$L\subset (M,\theta)$ is a two-part-exact Lagrangian cobordism.
By {\it a two-part-exact Lagrangian cobordism isotopy of $(L,\theta)$} we mean a smooth family
$\{ (L^\tau,\theta^\tau)\}_{0\leq\tau\leq T}$, so that $L^0=L$, $\theta^0=\theta$, and
for each $\tau\in [0,T]$,

\begin{itemize}

\item{} $(M,\theta^\tau)$ is an intergraphical exact symplectic cobordism based
on $(\Sigma,\xi)$,

\item{} $L^\tau\subset (M,\theta^\tau)$ is a two-part-exact Lagrangian cobordism.

\end{itemize}

The forms $d\theta^\tau$, $\tau\in [0,T]$, are also required to be identified with $d(c\theta_{can})$
for some $c>0$ (independent of $\tau$)
by a family of diffeomorphisms of $M$ depending smoothly on $\tau$ and equal to the identity on a common
collar neighborhood of the boundary of $M$.

We also say that this is a {\it two-part-exact Lagrangian cobordism isotopy from $(L,\theta)=(L^0,\theta^0)$ to
$(L^T,\theta^T)$}, or, if $\theta^\tau\equiv\theta$ is a constant family of forms,
{\it a two-part-exact Lagrangian cobordism isotopy from $L=L^0$ to
$L^T$}.
The function $\tau\mapsto \gap (\tau):= \gap_{L^\tau,\theta^\tau}$, $\tau\in [0,T]$,
is called {\it the gap function} of the isotopy.
The isotopy is said to be {\it fixed near the boundary} if all $(L^\tau,\theta^\tau)$ coincide
with $(L,\theta)$
on a collar neighborhood $U_\delta$ of $\partial M$ for a sufficiently small $\delta>0$
(see Definition~\ref{def-integraphical-exact-symplec-cob}).


\hfill


\subsection{Contactizations, non-degeneracy}
\label{subsec-prelims-contactization}

Following \cite{EES-TAMS} (see \cite{Gromov} for the original version of the notion,
cf. \cite{ALP}),
we say that a compatible almost complex structure $J$ on a (possibly non-compact) symplectic manifold
with a symplectic form $\omega$ has {\it finite geometry at infinity} if
the Riemannian metric $\omega (\cdot, J\cdot)$
admits a uniform positive
lower bound on the injectivity radius and a uniform upper bound on the sectional curvature.
A symplectic manifold is said to have {\it finite geometry at infinity} if it admits a compatible
almost complex structure that has finite geometry at infinity.

Further in the paper, we will focus on contactizations of exact symplectic manifolds that
have finite geometry at infinity -- i.e. on
contact manifolds of the form $\Sigma = P\times \R (z)$,
where $(P^{2n},d\vartheta)$, $n\in\Z_{> 0}$, is an exact symplectic manifold with a symplectic form $d\vartheta$
and finite geometry at infinity, and
where the contact structure $\xi$
on $\Sigma$ is defined by the contact form $\lambda_\vartheta = dz + \vartheta$. Note that the Reeb vector field of $\lambda_\vartheta$ is $\partial/\partial z$ and its flow has no periodic orbits.

We will say that a contact form $\lambda$
on such $(\Sigma,\xi)$ is {\it comfortable} if, for some $c>0$, it is obtained from
$c\lambda_\vartheta$ by a co-orientation-preserving contactomorphism of $(\Sigma,\xi)$.
In particular,
the Reeb flow of a comfortable contact form is defined for all times and has no periodic orbits (since the same is true
for the Reeb flow of $\lambda_\vartheta$).

Clearly, if $\lambda$ is comfortable, then so is $c\lambda$ for any $c>0$.

To the best of our knowledge, such $(\Sigma, \lambda)$ form the largest (up to a simultaneous diffeomorphism)
class of contact manifolds and contact forms
for which the details of the Legendrian contact homology theory are available in the published literature.

\hfill


\begin{exam}
\label{exam-1-jet-space}
{\rm
Let $Q$ be a closed manifold of positive dimension and let $P:= T^* Q$.
Let $\vartheta$ be equal to minus the canonical ``pdq" 1-form on $T^* Q$.

Then:

\begin{itemize}

\item{} $(P,d\vartheta)$ is an exact symplectic manifold with finite geometry at infinity (see \cite[Example 2.2]{EES-TAMS}).

\item{} $P\times \R (z) = J^1 Q$ is the 1-jet space of $Q$.

\item{} $\lambda_\vartheta: = dz + \vartheta$ is the standard contact form on $J^1 Q$ defining the standard contact structure
on $J^1 Q$.

\end{itemize}

}
\end{exam}


\hfill


\begin{exam}
\label{exam-compact-Liouville-domain}
{\rm
Recall that a {\it Liouville domain} is a compact exact symplectic manifold $(Y,d\vartheta)$ with boundary
such that the vector field $v$
on $Y$, defined by $d\vartheta (v,\cdot) = \vartheta (\cdot)$ and called a {\it Liouville vector field},
is non-zero and points outwards everywhere on $\partial Y$.

The form $\vartheta|_{\partial Y}$ is a contact form on $\partial Y$.
The flow of $v$ identifies a neighborhood of
$\partial Y$ in $Y$ with $\partial Y\times (1-\delta,1]$, for some small $\delta>0$, so that
$\vartheta$ is identified with $t(\vartheta|_{\partial Y})$ and $v$ is identified with $t\partial/\partial t$.
(Here $t$ is the coordinate on $(1-\delta,1]$).

Consequently, the triples $(Y,\vartheta,v)$ and
$\left(\partial Y \times [1,+\infty), t(\vartheta|_{\partial Y}), t\partial/\partial t\right)$
can be naturally glued along $\partial Y$ into a triple $(\bY,\bvartheta,\barv)$ of smooth objects.
The resulting exact symplectic manifold without boundary $(\bY, d\bvartheta)$ is called the {\it completion}
of $(Y,d\vartheta)$.

The completion $(\bY, d\bvartheta)$ has finite geometry at
infinity. Namely, one can take an almost complex structure $J$ on $\bY$ which is compatible
with $d\bvartheta$ and whose restriction to $\partial Y \times [1,+\infty)$ is invariant with respect
to the multiplicative $[1,+\infty)$-action of the second factor and which sends $t\partial/\partial t$
to the Reeb vector field of $\vartheta|_{\partial Y}$ everywhere on $\partial Y \times [1,+\infty)$.
A direct check shows that $J$ has finite geometry at infinity, and hence so does $(\bY,d\bvartheta)$.

}
\end{exam}


\hfill


An intergraphical exact symplectic cobordism $(\SSigma_{\mu^-,\mu^+},\theta)$, which is
based on $(\Sigma,\xi)$,
will be called
{\it comfortable}
if the forms $\mu^-,\mu^+$ are comfortable.
Clearly if $(\SSigma_{\mu^-,\mu^+},\theta)$ is comfortable, then so is $(\SSigma_{\mu^-,\mu^+},c\theta)$
for any $c>0$.

Let $\Lambda=\Lambda_0\sqcup\Lambda_1$ be a two-part Legendrian submanifold.

Denote the set of all Reeb chords of $\Lambda$ (with respect to $\lambda$) by
$\cR (\Lambda,\lambda)$.

We say that a Reeb chord of $\Lambda$ is an {\it $ij$-chord} for $i,j=0,1$ if it starts on $\Lambda_i$ and ends on $\Lambda_j$.
Denote the set of all $ij$-chords of $\Lambda$ by
$\cR_{ij} (\Lambda,\lambda)$.
The $ii$-chords will be called {\it pure}, while $ij$-chords for $i\neq j$ will be called {\it mixed}.

Note that there is a natural bijective correspondence between the sets
$\cR (\Lambda,\lambda)$ and $\cR (\Lambda,c\lambda)$ and, moreover, between the
sets $\cR_{ij} (\Lambda,\lambda)$ and $\cR_{ij} (\Lambda,c\lambda)$,
for any $i,j=0,1$ and any $c>0$.

We say that the pair $(\Lambda,\lambda)$ is {\it non-degenerate} if, with the Reeb flow of $\lambda$ being denoted
by $\{\psi_t\}$, the following conditions are satisfied:

\begin{itemize}

\item{} For each Reeb chord
$a: [0,T]\to \Sigma$, $a (t) = \psi_t
(a(0))$ (for all $t\in [0,T]$), $a(0), a(T)\in\Lambda$,
of $\Lambda$ with respect to $\lambda$,
the plane $d_{a(0)} \psi_T (T_{a(0)} \Lambda)$ is transverse to $T_{a(T)} \Lambda$ in the contact plane $\xi_{a(T)}$.

\item{} Each trajectory of the Reeb flow $\{ \psi_t (x)\}$, $-\infty < t < +\infty$, $x\in\Sigma$, intersects $\Lambda$
in at most two points.
    (Equivalently, the images of distinct Reeb chords are disjoint).
\end{itemize}

If $(\Lambda,\lambda)$ is non-degenerate, then

\begin{itemize}

\item{} $(\Lambda,c\lambda)$ is also non-degenerate for all $c>0$,

\item{} the set $\cR (\Lambda,\lambda)$ is finite.

\end{itemize}

A two-part-exact Lagrangian cobordism between $(\Lambda^\pm,\lambda^\pm)$ is said to be {\it non-degenerate}
if the pairs $(\Lambda^\pm,\lambda^\pm)$ are non-degenerate.

In case $\lambda=\lambda_\vartheta$, the non-degeneracy of a pair $(\Lambda,\lambda)$
can be described as follows.

Let
\[
\pr: \Sigma = P\times \R\to P
\]
be the natural projection.

Then $(\Lambda,\lambda)$ is non-degenerate if
\[
\pr|_\Lambda: \Lambda\to P
\]
is an immersion with only simple transverse double points. These double points are in bijective correspondence
with the Reeb chords of $(\Lambda,\lambda)$.

\bigskip
{\bf

Further on in the paper we will always assume that
$P$, $\vartheta$,
$\Sigma$, $\lambda_\vartheta$ and $\xi = \ker\lambda_\vartheta$ are as above and that
all contact forms on $\Sigma$ that we consider are comfortable.}
\bigskip


\hfill


\subsection{Persistence modules}
\label{subsec-persistence-modules}

In the paper we use some standard notions of the theory of persistence modules and a few new ones.
Let us spell out the precise setting that we use as well as the new notions that we introduce.
For a more detailed introduction to the theory of persistence modules
see e.g. \cite{Edelsbrunner-Harer}, \cite{Chazal-et-al-book}, \cite{Oudot-AMS2015-book}, or \cite{Polterovich-et-al-Persist-Homology-book}.

Let $\F=\Z_2$.


\hfill


\begin{defin}
\label{def-persistence-module}
{\rm
A {\it persistence module $(V,\pi)$ (or just $V$)} is given by
\[
\left(V=\{ V_t\}_{t\in (-\infty,+\infty]}, \pi = \{ \pi_{s,t}\}_{s,t\in (-\infty,+\infty], s\leq t}\right),
\]
where all $V_t$ are  $\F$-vector spaces, that are finite-dimensional for finite $t$ (but $V_\infty:=V_{+\infty}$ may be
infinite-dimensional)
and $\pi_{s,t} : V_s \to V_t$ are linear maps, so that

\smallskip
\noindent
(i) (Persistence) $\pi_{t,t} = Id$, $\pi_{s,r} = \pi_{t,r} \circ \pi_{s,t}$, for all $s, t, r \in (-\infty,+\infty]$, $s\leq t\leq r$.

\smallskip
\noindent
(ii) (Discrete spectrum and semicontinuity) There exists a (possibly empty) finite or countable discrete closed set of points
$\spec (V)\subset \R$,
called {\it the spectrum of $V$}, so that

\begin{itemize}

\item{}
For any $r \in \R$, $r\notin\spec (V)$, there exists a
neighborhood $U$ of $r$ in $\R$ such that $\pi_{s,t}$ is an isomorphism for all $s,t\in U$, $s \leq t$.

\item{}
For any $r\in \spec (V)$, there exists $\epsilon > 0$ such that $\pi_{s,t}$ is an isomorphism of vector spaces
for all $s,t\in [r -\epsilon,r]$, while for any $t\in (r,r +\epsilon]$, the map $\pi_{r,t}$ is
not an isomorphism.

\end{itemize}

\smallskip
\noindent
(iii) (Direct limit)
$V_\infty =\varinjlim\limits_{t\to +\infty} V_t$.

\smallskip
\noindent
(iv) (Semi-bounded support)
There exists $b\in\R$ such that $V_t=0$ for all $t\leq b$.


For $-\infty\leq a<+\infty$, we say that $V$ is {\it supported in $(a,+\infty)$} if
either $V_t = 0$ for all $t\in\R$
or
$\min \spec (V)\in (a,+\infty)$.

}
\end{defin}


\hfill


With this terminology, any persistence module is supported in $(-\infty,+\infty)$.


\hfill


\begin{defin}
\label{def-pers-modules-morphisms}
{\rm
Let $(V,\pi)$, $(\tV,\tpi)$ be persistence modules.
For $C\in\R$, a {\it morphism $\Phi: (V,\pi)\to (\tV,\tpi)$ of additive degree $C$}
is defined as a collection of linear maps
$\Phi_t: V_t\to \tV_{t-C}$, $t\in (-\infty,+\infty]$, such that
\[
\Phi_t\circ\pi_{s,t} = \tpi_{s-C,t-C}\circ \Phi_s\ \textrm{for all}\ s,t\in (-\infty,+\infty], s\leq t.
\]

We write $\Phi_\infty$ for $t=+\infty$.
}
\end{defin}


\hfill


With this terminology, the usual morphisms of persistence modules are the ones of additive degree $0$.

Clearly, a composition of morphisms of additive degrees $C$, $C'$ is a morphism of additive degree $C+C'$.

An invertible morphism (of additive degree $0$) is called {\it an isomorphism}.

We say that morphisms $\Phi: V\to W$, $\tPhi: \tV\to \tW$ (of the same additive degree) of persistence modules
are {\it RL-equivalent} if there exist isomorphisms $f: V\to\tV$, $g:W\to\tW$ such that $g\circ \Phi = \tPhi\circ f$.


\hfill


\begin{exam}
\label{exam-mult-shifted-morphism}
{\rm
Let $c>0$.

Let $V= (V,\pi)$ be a persistence module supported in $(0,+\infty)$.
Define a new persistence module $V^{[\times c]}$ supported in $(0,+\infty)$
by reparameterizing via $t\mapsto ct$ -- in particular,
\[
V^{[\times c]}_t:= V_{ct}.
\]

A morphism
\[
\Phi: V\to W
\]
of additive degree $C$ of persistence modules  supported in $(0,+\infty)$ induces a morphism
\[
\Phi^{[\times c]}: V^{[\times c]}\to W^{[\times c]}
\]
of additive degree $C/c$.

If $c\geq 1$, one can also define a morphism (of additive degree $0$)
\[
\Sh_V [\times c]: V\to V^{[\times c]},
\]
called {\it the multiplicative shift of $V$ by $c$}. It is formed by the maps
$\pi_{t,ct}: V_t\to V^{[\times c]}_t$ for $t\in (0,+\infty]$
and by the zero maps $V_t\to V^{[\times c]}_t$ for $t\in (-\infty,0]$.
}
\end{exam}


\hfill


\begin{exam}
\label{exam-add-shift}
{\rm
Let $V=(V,\pi)$ be a persistence module supported in $(0,+\infty)$ and let $C\geq 0$.

The collection of maps $\pi_{t,t+C}: V_t\to V_{t+C}$, $t\in (-\infty,+\infty]$, defines an endomorphism
\[
\Sh_V (+C): V\to V
\]
of additive degree $-C$. We will call it {\it an additive self-shift of $V$ by $C$}.

Sometimes we will use the notation $\Sh (+C)$ instead of $\Sh_V (+C)$.
}
\end{exam}


\hfill


\begin{rem}
\label{rem-add-shift-RL-equivalence}
{\rm
One easily checks that if
morphisms $\Phi: V\to W$ and $\tPhi: \tV\to \tW$ between persistence modules
are RL-equivalent (and, in particular, have the same additive degree),
then so are $\Sh_W (+C)\circ\Phi$ and $\Sh_\tW (+C)\circ\tPhi$ for any $C\geq 0$.

}
\end{rem}


\hfill


We will work with persistence modules supported in $(0,+\infty)$
and with the multiplicative scale for the indices
$t\in (0,+\infty)$ of the vector spaces in such persistence modules.
The identification with the usual scale is provided
by replacing all the indices $t\in \R$
for the usual persistence modules with $e^t$ for the persistence modules supported in $(0,+\infty)$.

Let $\J = (a_\J, b_\J]$, $-\infty< a_\J<b_\J\leq +\infty$. Define an {\it interval persistence module}
\[
(Q(\J), \pi) := \left(\{ Q(\J)\}_{t\in(-\infty,+\infty]} , \{\pi_{s,t}\}_{s,t\in(-\infty,+\infty],s\leq t}\right),
\]
as follows: $Q(\J)_t = \F$ for $t\in \J$ and $Q(\J)_t = 0$ for $t\in (-\infty,+\infty]\setminus\J$,
while the morphisms $\pi_{s,t}$ are the identity maps for $s,t\in\J$ and zeroes otherwise.

By the theory of persistence modules, each persistence module $V=\{V_t\}$ can be represented
as a locally finite (in $t$) direct sum of interval persistence modules.
Such a decomposition is called a {\it bar decomposition of $V$}; it may not be unique.
The resulting (possibly infinite, at most countable) collection of intervals is unique. It is called
{\it the barcode of $V$} and will be denoted by $\cB (V)$; the intervals in the collection
are called the {\it bars of $V$}, or the {\it bars in $\cB (V)$}.
In fact, the barcode $\cB (V)$ is a multiset: each bar may appear in the collection several, but only
finitely many, times.
Note that the barcode is locally finite -- that is, each point in $\R$ has an open neighborhood that
intersects at most finitely many bars. Persistence modules are isomorphic if and only if they have the same
barcode.

For a persistence module $V$ supported in $(0,+\infty)$, the left endpoints of
all the bars in $\cB (V)$ lie in $(0,+\infty)$ -- we will say that such a barcode is {\it supported in $(0,+\infty)$}.
We will perform the measurements
with such barcodes on the multiplicative (logarithmic) scale -- for instance, the
{\it multiplicative length} of a bar $\I=(a,b]$, $0<a<b$, is defined as
\[
|\I|^\times := b/a.
\]


\hfill


\begin{defin}
\label{def-lminc}
{\rm
For a barcode $\cB$ supported in $(0,+\infty)$ and for $c>1$, set
\[
\lminc (\cB) := \inf \{\ a_\I\ |\ \I = (a_\I,b_\I]\in \cB,\ |\I|^\times > c\ \}.
\]
Also, set
\[
\lmininfty (\cB) := \inf \{\ a_\I\ |\ \I = (a_\I,+\infty]\in \cB\ \}.
\]
If the infima above are taken over empty sets, we set them to be $+\infty$.
}
\end{defin}


\hfill


\begin{defin}
\label{def-birth-time}
{\rm

Given $\alpha\in V_t$, $t\in (0,+\infty]$, $\alpha\neq 0$, define {\it the birth-time of $\alpha$} as
\[
l_V (\alpha) := \inf\
\left\{\  s\in (0,+\infty],\ s\leq t\ |\ \alpha\in \Image \pi_{s,t}\ \right\}.
\]
Also set
\[
l_V (0):= 0.
\]
}
\end{defin}


\hfill


\begin{rem}
\label{rem-birth-time}
{\rm
$\ $

\smallskip
\noindent
1. For elements $\alpha\in V_\infty$, the birth-times $l_V (\alpha)$ are also known in symplectic topology as {\it spectral numbers}.

\bigskip
\noindent
2. If under a bar decomposition $V=\oplus_i V_i$ an element $\alpha$ is written as
$\alpha = \alpha_{i_1}+\ldots+\alpha_{i_k}$, for non-zero $\alpha_{i_1},\ldots,\alpha_{i_k}$,
then $l_V (\alpha)$ is the rightmost of the left ends
of the bars corresponding to the interval persistence modules $V_{i_1},\ldots, V_{i_k}$.

}
\end{rem}


\hfill


\subsection{The retreat of morphisms between persistence modules}
\label{subsec-retreat-morphisms-of-pers-modules}

Let $\Phi: V\to W$ be a morphism of additive degree $C_\Phi$ between persistence modules supported in $(0,+\infty)$.

We write $\alpha\in V$ for an element $\alpha\in V_t$ for some $t\in (-\infty,+\infty]$.


\hfill


\begin{defin}
\label{def-Upsilon}
{\rm
Define {\it the retreat} of $\Phi$ as
\[
\Ups (\Phi):=
\inf\ \left\{\ l_V (\alpha) - l_W (\Phi (\alpha))- C_\Phi\ |\ \alpha\in V,\ \Phi (\alpha)\neq 0\ \right\}.
\]
Similarly, define {\it the $\infty$-retreat} of $\Phi$
\[
\Ups_\infty (\Phi):=
\inf\ \left\{\ l_V (\alpha) - l_W (\Phi (\alpha))- C_\Phi\ |\ \alpha\in V_\infty,\ \Phi (\alpha)\neq 0\ \right\}.
\]
If the infima above are taken over empty sets, we set them to be $+\infty$.
}
\end{defin}


\hfill


In particular, $\Ups (\Phi)$ is finite if and only if the morphism $\Phi$ is non-trivial (i.e. not
equal identically to zero).

The proofs of the propositions below, listing the basic properties of the retreat,
are straightforward and are left as exercises for the reader.


\hfill


\begin{prop}
\label{prop-Ups-vs-Ups-infty}
\begin{equation}
\label{eqn-Ups-Phi-leq-Ups-infty-Phi}
\Ups (\Phi)\leq\Ups_\infty (\Phi)
\end{equation}
and
\begin{equation}
\label{eqn-Ups-infty-add-shift-composed-with-Phi-equals-Ups-infty-Phi}
\Ups_\infty \left(\Sh_W (+C)\circ\Phi\right) = \Ups_\infty (\Phi) + C,\ \textit{for all}\ C\geq 0.
\end{equation}
\end{prop}


\hfill


\begin{prop}
\label{prop-Ups-depends-only-on-isom-class-of-morphism}
If morphisms $\Phi: V\to W$, $\Psi: \tV\to\tW$ of persistence modules are RL-equi\-va\-lent,
then
\[
\Ups (\Phi) = \Ups (\Psi),
\]
\[
\Ups_\infty (\Phi) = \Ups_\infty (\Psi).
\]

\end{prop}


\hfill


\begin{prop}
\label{prop-retreat-mult-shift-of-a-morphism}
Let $c>0$.

Then
\[
\Ups \left( \Phi^{[\times c]}\right) = \frac{\Ups (\Phi)}{c}.
\]
\Qed
\end{prop}


\hfill


\begin{prop}
\label{prop-retreat-of-mult-shift-Sh}
Let
$c> 1$.

Then
\[
\Ups \left( \Sh_V [\times c]\right) = \frac{c-1}{c} \lminc \big(\cB(V)\big),
\]
\[
\Ups_\infty \left( \Sh_V [\times c]\right) = \frac{c-1}{c} \lmininfty \big(\cB(V)\big).
\]

If $c=1$ (i.e. $\Sh_V [\times c] = \Id$) and $V\neq 0$ (respectively, $V_\infty\neq 0$), then
$\Ups \left( \Sh_V [\times c]\right) = 0$ (respectively, $\Ups_\infty \left( \Sh_V [\times c]\right) = 0$).

\Qed
\end{prop}


\hfill


\begin{prop}
\label{prop-retreat-of-add-shift-Sh}
Let $C\geq 0$.

Then
\[
\Ups (\Sh_V (+C)) =
\begin{cases}
C, & \textrm{if } \Sh_V (+C)\neq 0,\\
+\infty, & \textrm{if } \Sh_V (+C)=0.
\end{cases}
\]
(By equating a morphism to $0$, we mean that it is identically equal to $0$).
\Qed
\end{prop}


\hfill


\begin{prop}
\label{prop-Upsilon-composition-of-morphisms}
$\ $

\smallskip
\noindent
A. If $\Phi,\Psi:V\to W$ are morphisms of the same additive degree, then
\[
\Ups(\Phi+\Psi)\geq \min\{\Ups(\Phi),\Ups(\Psi)\}.
\]

\smallskip
\noindent
B. For any morphisms $\Psi:U\to V$ and $\Phi:V\to W$, of arbitrary additive degrees, between persistence modules,
\[
\Ups(\Phi\Psi)\geq \Ups(\Phi)+\Ups(\Psi),
\]
\[
\Ups_\infty (\Phi\Psi)\geq \Ups_\infty (\Phi)+\Ups_\infty (\Psi).
\]
\Qed
\end{prop}


\hfill


The claim B in Proposition~\ref{prop-Upsilon-composition-of-morphisms} is crucial to our considerations. Its proof
is based on taking the infima in the left-hand side and in both terms in the right-hand side
of the following obvious identity:
\[
l_U \left(\alpha\right) - l_W \left(\Phi \Psi \left(\alpha\right)\right) - C_{\Phi \Psi} =
\]
\[
\left(l_U(\alpha) - l_V(\Psi(\alpha)) - C_\Psi\right)
+ \left(l_V(\Psi(\alpha)) - l_W(\Phi\Psi(\alpha)) - C_\Phi\right).
\]


\hfill


\section{Main results and applications}
\label{sec-lch-pers-modules}

Assume that $(\Sigma,\xi)$ is a contactization of an exact
symplectic manifold with finite geometry at infinity, as in Section~\ref{subsec-prelims-contactization}.

Let $\Lambda=\Lambda_0\sqcup\Lambda_1\subset (\Sigma,\xi)$ be a two-part
Legendrian submanifold and $\lambda$ a comfortable contact form on
$(\Sigma,\xi)$ so that the pair $(\Lambda,\lambda)$ is non-degenerate.

As we will see below in Section~\ref{subsec-Chekanov-Eliashberg-algebra-pers-mod}, one can use a version of Legendrian contact homology and an
additional choice
of an almost complex structure on $\Sigma\times\R_+$
to associate to the pair $(\Lambda,\lambda)$
a persistence module supported in $(0,+\infty)$.

The isomorphism class of this persistence module is independent of the choice of the almost complex structure
and will be denoted by $V(\Lambda_0,\Lambda_1,\lambda)$. Sometimes we will abuse the notation and denote by
$V(\Lambda_0,\Lambda_1,\lambda)$
the persistence module itself (that depends on an almost complex structure not appearing in the notation).
The corresponding barcode will be denoted by $\cB (\Lambda_0,\Lambda_1,\lambda)$.

The dimension $\dim V_\infty (\Lambda_0,\Lambda_1,\lambda)$ is well-defined and equals the number of
infinite bars in $\cB (\Lambda_0,\Lambda_1,\lambda)$ (where each bar is counted the number of times it appears in the
multiset $\cB (\Lambda_0,\Lambda_1,\lambda)$).
Moreover, $\dim V_\infty (\Lambda_0,\Lambda_1,\lambda)$ is preserved
by isotopies of $\Lambda=\Lambda_0\sqcup\Lambda_1$ in the class of two-part Legendrian submanifolds
(see Proposition~\ref{prop-dim-V-infty-indep-of-Leg-isotopy}).

Assume now that

\begin{itemize}

\item{} $(M,\theta)$ is a comfortable intergraphical exact symplectic cobordism, based on $(\Sigma,\xi)$,
between $(\Sigma^\pm,\lambda^\pm)$.

%

\item{} $L\subset (M,\theta)$ is a two-part-exact Lagrangian
cobordism between $(\Lambda^\pm,\lambda^\pm)$, where $\Lambda^\pm\subset (\Sigma,\xi)$ are
two-part Legendrian submanifolds.

\item{} the pairs $(\Lambda^\pm,\lambda^\pm)$
 are non-degenerate (and thus, with our terminology,
$L$ is non-degenerate).

\end{itemize}

The same version of Legendrian contact homology as above and an additional choice
of an almost complex structure on $\Sigma\times\R_+$ associate to the pair $(L,\theta)$
a morphism $V(\Lambda^+_0,\Lambda^+_1,\lambda^+)\to V(\Lambda^-_0,\Lambda^-_1,\lambda^-)$ of additive degree $\gap_{L,\theta}$
between the persistence modules. The morphism depends on the choice of the almost complex structure,
but its RL-equivalence class does not --
see
Proposition~\ref{prop-Phi-L-I-indep-of-I-up-to-equivalence}. We will denote the latter RL-equivalence class by $\Phi^{L,\theta}_*$. Sometimes
we will abuse the notation and denote by $\Phi^{L,\theta}_*$ the morphism itself
(that depends on an almost complex structure not appearing in the notation).

Moreover, by Remark~\ref{rem-add-shift-RL-equivalence},
for any $C\geq 0$, the RL-equivalence class of the composition of the latter morphism
with the additive shift by $C$ of its target is also independent of the choice of the almost complex
structure; we will denote it by
\[
\Sh (+C)\circ \Phi^{L,\theta}_*.
\]
By Remark~\ref{rem-add-shift-RL-equivalence}
and Proposition~\ref{prop-Ups-depends-only-on-isom-class-of-morphism}, the retreat
\[
\Ups \left( \Sh (+C)\circ \Phi^{L,\theta}_*\right)
\]
is well-defined.


\hfill


\subsection{The results on the retreat of two-part-exact Lag\-ran\-gian cobordisms}
\label{subsec-main-results-on-retreat}

Here are our main results concerning the retreat $\Ups \left( \Sh (+C)\circ \Phi^{L,\theta}_*\right)$ defined above.
(The order of these statements may not match the order of their proofs).


\hfill


\begin{thm}
\label{thm-Ups-of-cobordism-positive}
With the setting as above,

\bigskip
\noindent
A. $\Ups \left(\Phi^{L,\theta}_*\right) >0$.

\bigskip
\noindent
B. $\Ups \left(\Phi^{L,c\theta}_*\right) = c \Ups \left(\Phi^{L,\theta}_*\right)$ for any $c>0$.

\end{thm}


\hfill


\begin{thm}
\label{thm-trivial-morphism}
Let $0<s_- <s_+$, $c:= s_+/s_-$.

Assume that in the setting above

\begin{itemize}


\item{}
$\lambda$ is a comfortable contact form on $(\Sigma,\xi)$
and $(M=\Sigma\times [s_-,s_+], \theta=s\lambda)$
is a trivial intergraphical exact symplectic cobordism,

\item{} $\Lambda^\pm = \Lambda$,

\item{} $L = \Lambda\times [s_-,s_+]\subset (\Sigma\times [s_-,s_+], s\lambda)$ is a
non-degenerate
trivial two-part-exact
Lagrangian cobordism.

\end{itemize}

Then
\[
\Ups \left( \Phi^{L,s\lambda}_*\right) = (s_+ - s_-) \lminc \left(\cB (\Lambda_0,\Lambda_1,\lambda)\right),
\]
\[
\Ups_\infty \left( \Phi^{L,s\lambda}_*\right) = (s_+ - s_-) \lmininfty \left(\cB (\Lambda_0,\Lambda_1,\lambda)\right).
\]

\end{thm}


\hfill


The following theorem is central to our considerations.


\hfill


\begin{thm}
\label{thm-two-part-exact-Lagr-isotopy-morphisms-of-pmodules}
With the setting as above, assume that
$\{ (L^\tau,\theta^\tau)\}_{0\leq\tau\leq T}$
is a two-part-exact Lagrangian isotopy in $M$
fixed near the boundary, and that $L^0$ and $L^T$ are non-degenerate.

Let $\gap (\tau) := \gap_{L^\tau,\theta^\tau}$ be the gap function of the isotopy.
Set:
\[
\gap_{min} := \min_{\tau\in [0,T]} \gap (\tau),
\]
\[
C_0 := \gap (0) - \gap_{min},\ C_T := \gap (T) - \gap_{min}.
\]

Then
\[
\Ups_\infty \left(\Phi^{L^0,\theta^0}_*\right) + C_0\geq
\Ups \left(\Sh (+C_0)\circ \Phi^{L^0,\theta^0}_*\right) \geq
\Ups \left(\Phi^{L^T,\theta^T}_* \right) + C_T,
\]
\[
\Ups_\infty \left(\Phi^{L^0,\theta^0}_*\right) + C_0 =
\Ups_\infty \left(\Phi^{L^T,\theta^T}_*\right) + C_T.
\]

\end{thm}


\hfill


The following result is also crucial for our considerations: it
describes the behavior of $\Ups \left(\Phi^{L,\theta}_*\right)$,
$\Ups_\infty \left(\Phi^{L,\theta}_*\right)$ under concatenations.


\hfill


\begin{thm}
\label{thm-Ups-concatenation-of-cobordisms}
Assume that
$(M,\theta)$, $(\tM,\ttheta)$
are comfortable intergraphical exact symplectic cobordisms based on $(\Sigma,\xi)$
and
$L\subset (M,\theta)$, $\tL\subset (\tM,\ttheta)$
are non-degenerate two-part-exact Lagrangian cobordisms
such that the concatenation $L\sharp\tL\subset (M\sharp \tM, \theta\sharp \ttheta)$
is well-defined.

Then
\[
\Ups \left(\Phi^{L\sharp\tL,\theta\sharp\ttheta}_*\right) \geq \Ups\left(\Phi^{L,\theta}_*\right)
+ \Ups\left( \Phi^{\tL,\ttheta}_*\right),
\]
\[
\Ups_\infty \left(\Phi^{L\sharp\tL,\theta\sharp\ttheta}_*\right) \geq \Ups_\infty \left(\Phi^{L,\theta}_*\right)
+ \Ups_\infty \left( \Phi^{\tL,\ttheta}_*\right).
\]

\end{thm}


\hfill


\subsection{Application: a restriction on the existence of two-part-exact Lagrangian cobordism isotopies}
\label{subsec-appl-restriction-on-gap-fcns}

Here is a corollary of previous results providing an obstruction to the existence
of two-part-exact Lagrangian isotopies, fixed near the boundary, with certain kinds of gap functions.

Let us introduce the following definitions in the setting as above.


\hfill


\begin{defin}
\label{def-gap-raising-isotopy}
{\rm
By a {\it gap-raising Lagrangian isotopy} we mean a
two-part-exact Lagrangian cobordism isotopy $\{ L^\tau,\theta^\tau\}_{0\leq\tau\leq T}$,
fixed near the boundary, with a gap function $\gap (\tau)$ that satisfies the following condition:
\[
\gap (0) = \min_{\tau\in [0,T]} \gap (\tau).
\]
}
\end{defin}


\hfill


The following result describes a general phenomenon whose sample application was
presented in Theorem~\ref{thm-intro-application}.


\hfill


\begin{thm}
\label{thm-obstruction-to-Lagr-isotopies}
Assume that:

\begin{itemize}

\item{}
$L,\tL\subset (\Sigma\times [s_-,s_+],s\lambda)$ are non-degenerate two-part-exact Lagrangian cobordisms
such that there exists a two-part-exact Lagrangian isotopy from $L$ to $\tL$, fixed near the boundary.

\item{}
$\tLambda = \tLambda_0\sqcup \tLambda_1\subset (\Sigma,\xi)$
is a non-degenerate two-part Legendrian submanifold.

\item{} $\tLambda_i$, $L_i$ and $\tL_i$, $i=0,1$,  are all connected.

\item{}
$0< s_- < s_0 < s_1 < s_+$.

\item{} $\tc:= s_1/s_0$.

\item{} $\lmintc \big(\cB (\tLambda_0,\tLambda_1,\lambda)\big) < +\infty$ (i.e. $\cB (\tLambda_0,\tLambda_1,\lambda)$
contains bars of multiplicative length greater than $\tc$).

\item{} The intersection of $\tL$
with $\Sigma\times [s_0-\delta,s_1+\delta]$, for some $\delta>0$,
is the non-degenerate trivial two-part-exact Lagrangian cobordism
$\tLambda\times [s_0-\delta,s_1+\delta]\subset (\Sigma\times [s_-,s_+],s\lambda)$.

\end{itemize}

Then
\[
\gap_{\tL,s\lambda} <
\gap_{L,s\lambda} +
\Ups_\infty\left(\Phi^{L,s\lambda}_*\right) -
(s_1-s_0) \lmintc \big(\cB (\tLambda_0,\tLambda_1,\lambda)\big),
\]
and if the Lagrangian isotopy from $L$ to $\tL$ is gap-raising, then
\[
\gap_{\tL,s\lambda} <
\gap_{L,s\lambda} +
\Ups\left(\Phi^{L,s\lambda}_*\right) -
(s_1-s_0) \lmintc \big(\cB (\tLambda_0,\tLambda_1,\lambda)\big).
\]
\end{thm}


\hfill


\begin{cor}
\label{cor-obstruction-to-Lagr-isotopies-trivial-cobordism}
In the setting and under the assumptions of Theorem~\ref{thm-obstruction-to-Lagr-isotopies},
assume in addition that

\begin{itemize}

\item{}
$\Lambda = \Lambda_0\sqcup \Lambda_1\subset (\Sigma,\xi)$
is a non-degenerate two-part Legendrian submanifold,

\item{}
$L:= \Lambda\times [s_-,s_+]\subset (\Sigma\times [s_-,s_+],s\lambda)$,
meaning that $L$ is a non-degenerate trivial two-part-exact Lagrangian cobordism,

\item{} $c:= s_+/s_-$,

\item{} $\lminc \big(\cB (\Lambda_0,\Lambda_1,\lambda)\big)<+\infty$
(i.e. $\cB (\Lambda_0,\Lambda_1,\lambda)$
contains bars of multiplicative length greater than $c$).

\end{itemize}

Then
\[
\gap_{\tL,s\lambda} <
(s_+-s_-) \lmininfty \big(\cB (\Lambda_0,\Lambda_1,\lambda)\big) -
(s_1-s_0) \lmintc \big(\cB (\tLambda_0,\tLambda_1,\lambda)\big),
\]
and if the Lagrangian isotopy from $L$ to $\tL$ is gap-raising, then
\[
0\leq \gap_{\tL,s\lambda} <
(s_+-s_-) \lminc \big(\cB (\Lambda_0,\Lambda_1,\lambda)\big) -
(s_1-s_0) \lmintc \big(\cB (\tLambda_0,\tLambda_1,\lambda)\big).
\]
\end{cor}

\begin{proof}
For the trivial $L$, we have
\[
\gap_{L,s\lambda} = 0,
\]
and if the Lagrangian isotopy from $L$ to $\tL$ is gap-raising, then also
\[
\gap_{\tL,s\lambda} \geq 0.
\]
Also, by Theorem~\ref{thm-trivial-morphism},
\[
\Ups_\infty\left(\Phi^{L,s\lambda}_*\right) = (s_+-s_-) \lmininfty \big(\cB (\Lambda_0,\Lambda_1,\lambda)\big),
\]
\[
\Ups \left(\Phi^{L,s\lambda}_*\right) = (s_+-s_-) \lminc \big(\cB (\Lambda_0,\Lambda_1,\lambda)\big).
\]
Now the corollary follows directly from Theorem~\ref{thm-obstruction-to-Lagr-isotopies}.
\end{proof}


\hfill


Figure~\ref{fig-cobordism-schematic-intro}
illustrates a special case of the setting of
Corollary~\ref{cor-obstruction-to-Lagr-isotopies-trivial-cobordism}.

%
%
%
%
%
%
%
%
%
%
%
%
%

Theorem~\ref{thm-intro-application} in Section~\ref{sec-intro-and-sample-application} will be deduced from
Corollary~\ref{cor-obstruction-to-Lagr-isotopies-trivial-cobordism} and from the computation of
$\lminc$, $\lmintc$ and $\lmininfty$ for the corresponding two-part Legendrian submanifolds.


\hfill


\begin{rem}
\label{rem-examples-of-barcodes}
{\rm
For the two-part Legendrian submanifolds in Theorem~\ref{thm-intro-application} the barcodes
will contain infinite bars and, moreover, the smallest left end of the bars in each of these barcodes will
be the left end of an infinite bar, meaning that for these barcodes
$\lminc=\lmintc=\lmininfty$.

Let us note that there do exist examples of two-part Legendrian submanifolds in the contactizations of exact symplectic
manifolds with finite geometry at infinity, whose barcodes contain only finite bars, meaning that
$\lmininfty = +\infty$. For such two-part Legendrian submanifolds, the inequalities with $\Ups_\infty$ and $\lmininfty$
in Theorem~\ref{thm-obstruction-to-Lagr-isotopies} and Corollary~\ref{cor-obstruction-to-Lagr-isotopies-trivial-cobordism}
are trivial but the inequalities with $\Ups$ and $\lminc$, $\lmintc$ for finite $c,\tc$ do provide non-trivial obstructions.

Let us briefly describe examples of such two-part Legendrian submanifolds.

First, let $(\Sigma,\xi)$ be the standard contact $\R^3 (p,q,z)$ with the standard contact form $\lambda= dz-pdq$
on it. It is the contactization of the exact symplectic manifold $(P=\R^2,d\vartheta)$, where $\vartheta=-pdq$, with finite geometry at
infinity.
Consider the Lagrangian projection $pr: \R^3 (p,q,z)\to \R^2 (p,q)$. Consider
a two-part Legendrian submanifold $\Lambda_0\sqcup\Lambda_1\subset \R^3$ that represents a topological unlink
and whose Lagrangian diagram $pr (\Lambda_0)\sqcup pr (\Lambda_1)$ is as in Figure~\ref{fig-leg-unlink}.
The part $\Lambda_1$ lies above the part $\Lambda_0$ with respect to the $z$-coordinate.

\begin{figure}[htbp]
\centering
\includegraphics[width=0.7\textwidth]{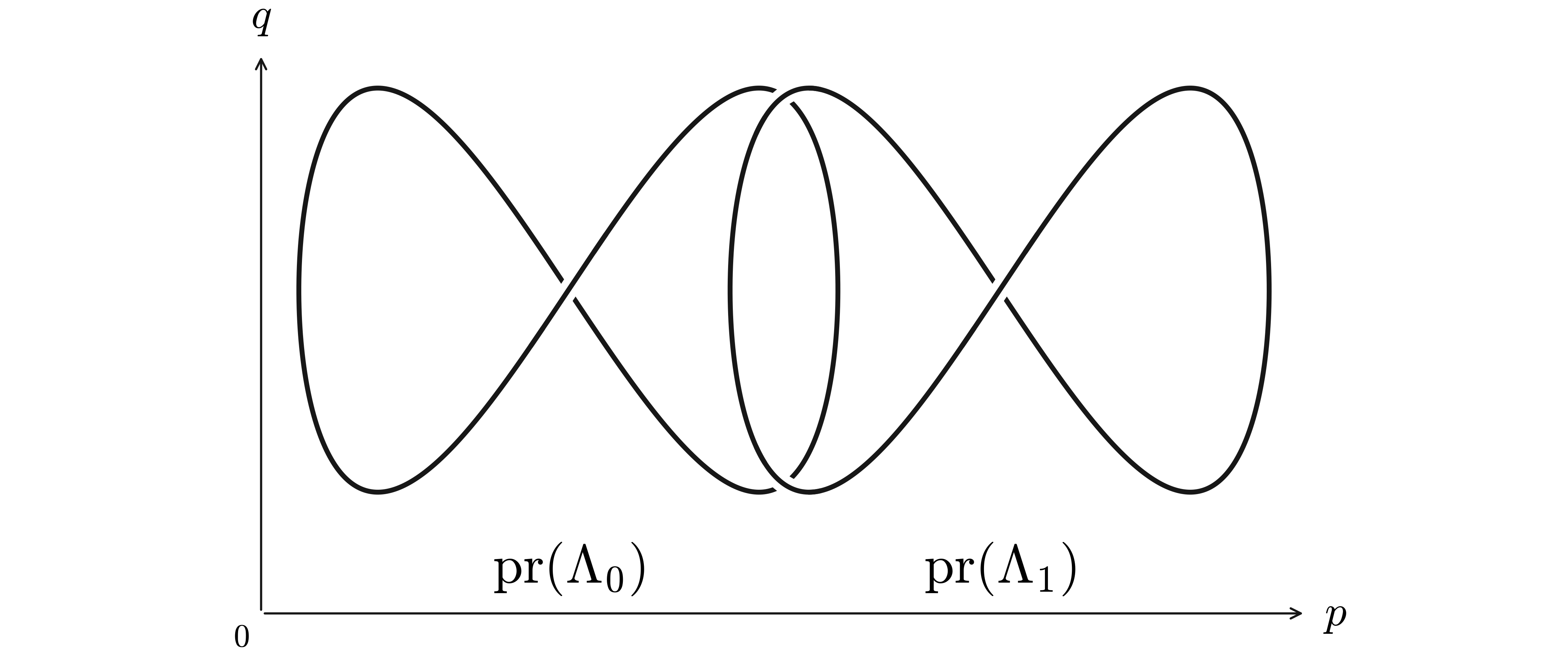}
\caption{The Lagrangian projection of a Legendrian unlink $\Lambda_0\sqcup\Lambda_1$}
\label{fig-leg-unlink}
\end{figure}

The pair $(\Lambda,\lambda)$ is non-degenerate and there are exactly 4 Reeb chords of $\Lambda$ with respect to $\lambda$:
one $00$-chord, one $11$-chord and two $01$-chords. One can show that the barcode $\cB (\Lambda_0,\Lambda_1,\lambda)$
consists of infinitely many finite bars -- see \cite{Avinoam-MSc-thesis}.

Let us now outline another example in which the barcode consists of {\sl exactly one finite bar}.

Consider the closed unit cotangent disk bundle $D^* \SP^1$ of $\SP^1$. It can be viewed as a compact Liouville
domain (see Example~\ref{exam-compact-Liouville-domain} for the definition).
Consider a boundary connected sum $Y$ of two copies of $D^* \SP^1$ obtained by attaching a Weinstein $1$-handle
between their boundaries; this is again a compact Liouville domain -- see \cite{Eliashberg-Topological-characterization-of-Stein-manifolds},
\cite{Weinstein-Hokkaido}. Its completion, to be denoted by $(P,d\vartheta)$,
is an exact symplectic manifold with finite geometry
at infinity (see Example~\ref{exam-compact-Liouville-domain}). In geometric terms, $P$ is $\R^2$ with two punctures
equipped with a special (non-standard) exact symplectic form.

Set $\Sigma:=P\times\R$, $\lambda := dz + \vartheta$.

The two zero-sections in $D^* \SP^1$ define disjoint embedded
exact circles
\[
\Gamma_0,\Gamma'_1\subset P.
\]
Apply a compactly supported Hamiltonian isotopy $\{\phi_t\}$ supported in the interior of $Y$ and
move $\Gamma'_1$ into an embedded circle $\Gamma_1$ that has two intersection
points $x,y$ with $\Gamma_0$, so that these points form
a cancelling pair of transverse intersections with
$\Gamma_0$. The arcs of the circles $\Gamma_0$, $\Gamma_1$
bound an embedded bigon $B$, with the vertices $x,y$, see Figure~\ref{fig-one-finite-bar}.

\begin{figure}[htbp]
\centering
\includegraphics[width=0.7\textwidth]{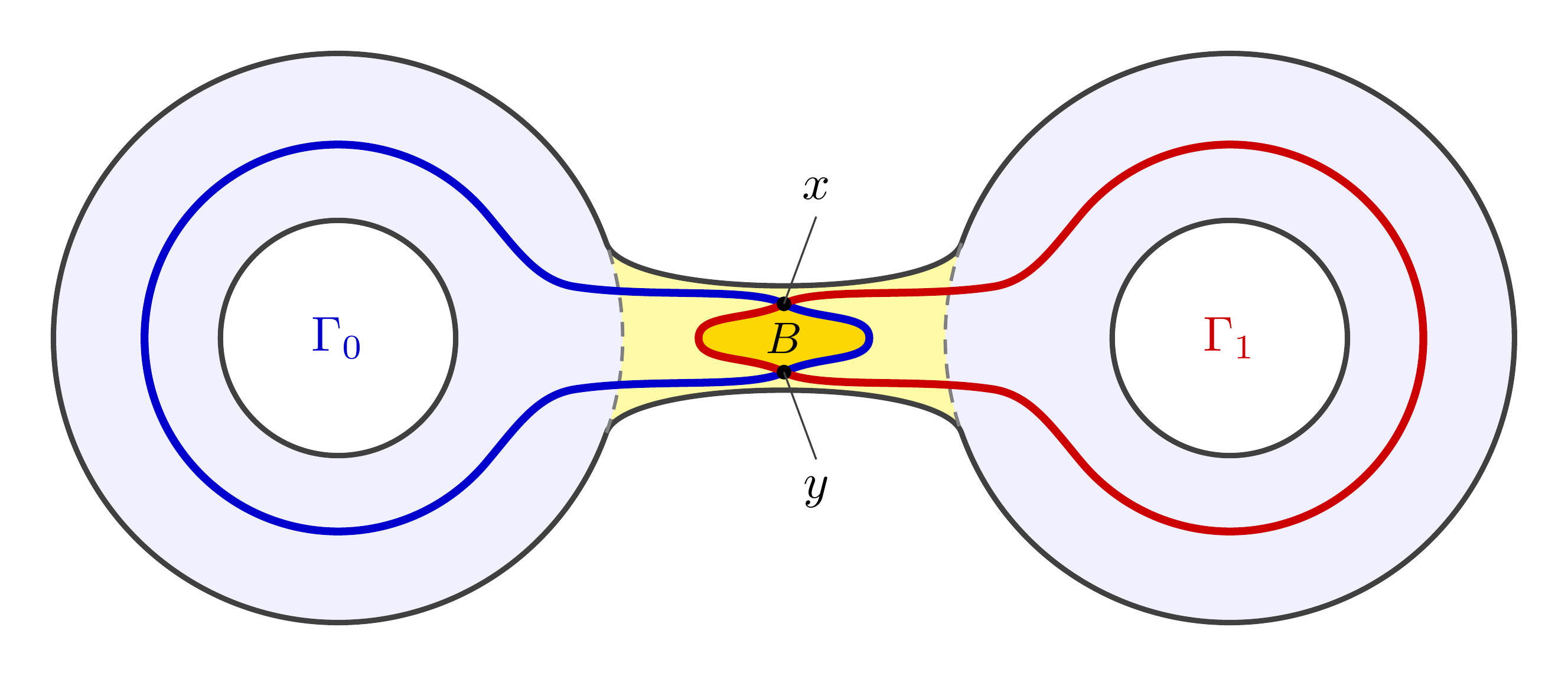}
\caption{A pair of Legendrians with exactly one finite bar}
\label{fig-one-finite-bar}
\end{figure}


Choose primitives $f_i: \Gamma_i\to\R$, $i=0,1$, such that $\vartheta|_{\Gamma_i}=df_i$.
For $a_i\in\R$, $i=0,1$,
define Legendrian lifts $\Lambda_0,\Lambda_1\subset (\Sigma, \ker\lambda)$ of $\Gamma_i$ by
\[
 \Lambda_i := \{(u,a_i-f_i(u)): u\in\Gamma_i\},\  i=0,1.
\]
If $a_1-a_0$ is chosen to be a sufficiently large positive number, then, as one can easily check:
\begin{itemize}

\item{}
$\Lambda:=\Lambda_0\sqcup\Lambda_1$ is a two-part Legendrian submanifold of $(\Sigma, \ker\lambda)$,

\item{} the pair $(\Lambda,\lambda)$ is non-degenerate,

\item{} there are exactly 2 Reeb chords of $\Lambda$ with respect to $\lambda$; both of them
are $01$-chords.
\end{itemize}

There is a Legendrian isotopy that moves $\Lambda$ into a two-part Legendrian submanifold whose Lagrangian projection
is the union of the disjoint curves $\Gamma_0$ and $\Gamma'_1$ and which therefore
has no Reeb chords at all (since $\Gamma_0$ and $\Gamma'_1$ are embedded and disjoint), and, in particular,
no $01$-chords. Recall that $\dim V_\infty (\Lambda_0,\Lambda_1,\lambda)$ is preserved
by isotopies of $\Lambda=\Lambda_0\sqcup\Lambda_1$ in the class of two-part Legendrian submanifolds
(see Proposition~\ref{prop-dim-V-infty-indep-of-Leg-isotopy}). Consequently,
$\dim V_\infty (\Lambda_0,\Lambda_1,\lambda) = 0$, meaning that
the barcode $\cB (\Lambda_0,\Lambda_1,\lambda)$ has no infinite bars.
Since there are exactly 2 Reeb chords of $\Lambda$ with respect to $\lambda$, and both of them
are $01$-chords, one easily gets that $\cB (\Lambda_0,\Lambda_1,\lambda)$ contains exactly one finite bar whose additive length equals the area of the bigon $B$.
}
\end{rem}


\hfill


\subsection{Applications to contact and symplectic dynamics}
\label{subsec-appls-to-dynamics}

Let us now outline an application of the retreat $\Ups\left(\Phi^{L,\theta}_*\right)$ to contact and symplectic Hamiltonian
dynamics.

Let us introduce the following definitions in the setting as above.


\hfill


\begin{defin}
\label{def-df}
{\rm
For a two-part-exact Lagrangian cobordism $L=L_0\sqcup L_1\subset (M,\theta)$, write:
\[
\df (L_0,L_1,\theta):= \sup_{(\tL,\ttheta)} (\gap_{\tL,\ttheta} - \gap_{L,\theta}),
\]
where the supremum is taken over all cobordism forms $\ttheta$ on $M$ and all
two-part-exact Lagrangian cobordisms
$\tL\subset (M,\ttheta)$ such that there exists a gap-raising Lagrangian isotopy from
$(L,\theta)$ to $(\tL,\ttheta)$.
}
\end{defin}


\hfill


A simple case in which $\df (L_0,L_1,\theta)$ can be computed is the
setting of Example~\ref{exam-Sigma-is-R-intro}, where $\dim \Sigma = 1$. In that case one easily sees
by elementary area considerations that
$\df (L_0,L_1,\theta) = (s_+ - s_-) (z_2 - z_1)$, and in particular $\df (L_0,L_1,\theta)$ is finite.

In higher dimensions, the finiteness of $\df (L_0,L_1,\theta)$ is a symplectic rigidity phenomenon. Our results
show that in order to prove the finiteness of $\df (L_0,L_1,\theta)$ it is enough to verify that the morphism
$\Phi^{L,\theta}_*$ of persistence modules is non-trivial.

Indeed, for a gap-raising Lagrangian isotopy,
we have $C_0 = \gap (0)- \gap (0) = 0$. Consequently, Theorem~\ref{thm-two-part-exact-Lagr-isotopy-morphisms-of-pmodules}
yields the following inequality:
\[
\df (L_0,L_1,\theta)\leq \Ups\left(\Phi^{L,\theta}_*\right),
\]
where the right-hand side is finite if and only if the morphism
$\Phi^{L,\theta}_*$ is non-trivial.

The finiteness of $\df (L_0,L_1,\theta)$ has important implications for contact and symplectic Hamiltonian dynamics.
Namely, one can show (see \cite{EP-big})
that if $\df (L_0,L_1,\theta) < +\infty$, then the quadruple $(L_0,L_1,\Sigma^-, \Sigma^+)$
has strong dynamical properties: roughly speaking, the symplectic Hamiltonian
flows of appropriate Hamiltonian functions whose level sets
separate the sets $\Sigma^-$ and $\Sigma^+$ (respectively, $L_0$ and $L_1$) have
trajectories connecting the other two sets in the quadruple. The result actually
provides the direction of the trajectory and an explicit upper bound on its time-length.
The proof is based on Poisson bracket invariants and is very similar
to the proof of Theorem 4.9 in \cite{EP-JFPTA2022} (where only trivial
Lagrangian cobordisms were considered). Let us outline it briefly.

Equip $M$ with the symplectic form $d\theta$. For a closed set
$X\subset M$, let $\Op (X)$ denote some open neighborhood of $X$. If $X=\Sigma^\pm$, we assume that
$\Op (X)$ contains a collar neighborhood of the form $\Sigma^+\times [1-\delta,1]$, respectively $\Sigma^-\times [1,1+\delta]$,
for some sufficiently small $\delta>0$.
Consider smooth functions $F,G: M\to [0,1]$ such that
\begin{itemize}

\item{} $F|_{\Op (L_0)} \equiv 0$, $F|_{\Op (L_1)} \equiv 1$,
$G|_{\Op (\Sigma^-)} \equiv 0$, $G|_{\Op (\Sigma^+)} \equiv 1$,

\item{} $\sup_M \{F,G\} < +\infty$ and the vector field $F\sgrad G$ on $M$ is complete.

\end{itemize}

Assume first that $\sup_M \{F,G\} >0$.
Pick any $0<T < 1/\sup_M \{ F,G\}$ and
consider a family
\[
\theta^\tau := \theta + \tau F dG, \ \tau\in [0,T],
\]
of 1-forms on $M$.
One can check that $d\theta^\tau$ is symplectic for all $\tau\in [0,T]$,
and $\{ L^\tau\equiv L, \theta^\tau\}$, $\tau\in [0,T]$, is a gap-raising
Lagrangian isotopy. Consequently,
\[
1/\sup_M \{ F,G\}\leq \df (L_0,L_1,\theta)\leq \Ups\left(\Phi^{L,\theta}_*\right),
\]
and therefore (since, by part A of Theorem~\ref{thm-Ups-of-cobordism-positive}, $\Ups\left(\Phi^{L,\theta}_*\right)>0$)
\[
\sup_M \{ F,G\}\geq 1/\Ups\left(\Phi^{L,\theta}_*\right).
\]
Assume now that $\sup_M \{F,G\} = 0$. Then the construction above yields
\[
\df (L_0,L_1,\theta) = \Ups\left(\Phi^{L,\theta}_*\right) = +\infty.
\]

Since this holds for all $F,G$ as above, we obtain that the
Poisson bracket invariants\footnote{By the Poisson bracket invariants we mean the ones
that are denoted in \cite{EP-JFPTA2022} by $pb^+_M$ and $pb^+_{M,comp}$ for
$M$ equipped with the symplectic form $d\theta$.} of the quadruple $(L_0,L_1,\Sigma^-, \Sigma^+)$
are positive and bounded from below by the number $1/\Ups\left(\Phi^{L,\theta}_*\right)$, which is positive, as long as
$\Ups\left(\Phi^{L,\theta}_*\right) < +\infty$.

Then, by
Proposition 2.6 in \cite{EP-JFPTA2022}, the positivity of the Poisson bracket invariants
yields the existence of the connecting trajectories along with an upper bound on their time-length.
The existence of connecting trajectories from $L_0$ to $L_1$
for trivial two-part-exact Lagrangian cobordisms $L=L_0\sqcup L_1 = (\Lambda_0\sqcup\Lambda_1)\times [s_-,s_+]$
yields the existence of trajectories from $\Lambda_0$ to $\Lambda_1$ for contact Hamiltonian flows on
$\Sigma$ -- see \cite{EP-JFPTA2022}.


\hfill


\subsection{Additional comments}
\label{subsec-comments-on-proofs}

The following comments on the results above are in order.


\hfill


\begin{rem}
\label{rem-deg-case-upper-semicontinuity}
{\rm
One can show that the functions
\[
\Lambda=\Lambda_0\sqcup\Lambda_1\mapsto \lminc \left( \cB \left( \Lambda_0,\Lambda_1,\lambda\right)\right)
\]
and
\[
L\mapsto \Ups \left(\Phi^{L,\theta}_*\right),
\]
for non-degenerate $\Lambda$, $L$,
are upper semicontinuous with respect to the $C^1$-topology on the domain. Degenerate Legendrian
submanifolds and Lag\-ran\-gian cobordisms
can be approximated in the $C^1$-topology by non-degenerate ones. Consequently, the functions above
extend as (possibly $+\infty$-valued) upper semicontinuous functions on the space of all,
possibly degenerate, $\Lambda$ and $L$. One can show that the results in
Sections~\ref{subsec-appl-restriction-on-gap-fcns}
and \ref{subsec-appls-to-dynamics}
extend to degenerate two-part-exact Lagrangian cobordisms (with strict inequalities being replaced by
non-strict ones).

For details see \cite{EP-big}.
}
\end{rem}

%
%
%

%
%
%
%


\hfill


\begin{rem}
\label{rem-sabloff-traynor}
{\rm
Let us compare our results with the results of Sabloff and Traynor in \cite{Sabloff-Traynor}.

\bigskip
\noindent
1. As we have already mentioned above, in our setting, a Lagrangian cobordism
is a compact submanifold with boundary
(of an intergraphical exact symplectic cobordism), while in \cite{Sabloff-Traynor}
a Lagrangian cobordism is a submanifold
{\sl without boundary} (of the symplectization of a 1-jet space of a manifold).
In other words, the Sabloff-Traynor Lagrangian cobordisms are the completions of our
Lagrangian cobordisms.

An (exact) Lagrangian isotopy of an (exact) Lagrangian cobordism in
\cite{Sabloff-Traynor} lives in the whole symplectization and may have arbitrarily large
(but compact) support, while in our case the isotopies of Lagrangian cobordisms
live in a given exact symplectic cobordism
with fixed boundaries and have to be supported away from these boundaries.

\bigskip
\noindent
2. In Section~\ref{subsec-appl-restriction-on-gap-fcns}, we study whether a given (two-part-exact) Lagrangian cobordism
can be isotoped (among two-part-exact Lagrangian cobordisms with the same boundary) into a Lagrangian cobordism
which is trivial in a {\sl given} exact symplectic subcobordism $\Sigma\times [s_0,s_1]$ (an ``$[s_0,s_1]$-layer")
in a given trivial exact symplectic cobordism $\Sigma\times [s_-,s_+]$.

The question that Sabloff and Traynor study in \cite{Sabloff-Traynor} can be described in our terms as  the question
about the smallest $\log (s_1/s_0)$ such that a given (exact) Lagrangian cobordism
can be isotoped (among exact Lagrangian cobordisms with the cylindrical ends) to a Lagrangian cobordism
which is trivial {\sl outside some} ``$[s_0,s_1]$-layer".

\bigskip
\noindent
3. The methods used in \cite{Sabloff-Traynor} to bound $\log (s_1/s_0)$ from below use a different version
of the Legendrian contact homology (namely, the linearized
Legendrian contact homology) filtered by the action of Reeb chords and a quantity analogous to what we call
$\Ups_\infty$ (a version of the retreat). In our results, we use not only $\Ups_\infty$ but $\Ups$ itself
which allows us to use the information ``contained in finite bars".

}
\end{rem}

%
%
%
%
%
%


\hfill


\section{Legendrian contact homology}
\label{sec-lch}

In this section we discuss the Legendrian contact homology constructions, mentioned briefly
in Section~\ref{sec-lch-pers-modules}.

We assume that $(\Sigma,\xi)$ is a contactization of an exact
symplectic manifold with finite geometry at infinity, as in Section~\ref{subsec-prelims-contactization}.


\hfill


\subsection{Chekanov-Eliashberg algebra and the corresponding persistence modules}
\label{subsec-Chekanov-Eliashberg-algebra-pers-mod}


Let $\lambda$ be a comfortable
contact form on $(\Sigma,\xi)$.
Assume $\Lambda\subset (\Sigma,\xi)$ is a Legendrian submanifold so that the pair $(\Lambda,\lambda)$ is non-degenerate.
Let us identify $(S\Sigma,\theta_{can})$ with $(\Sigma\times\R_+,s\lambda)$ as in
Section~\ref{subsec-prelims-intergraphical}. Under this identification, the cylinder $S\Lambda$ is identified with $\Lambda\times\R_+$.

Let us outline the definition of the Chekanov-Eliashberg differential algebra associated to
$(\Lambda,\lambda)$ --
see \cite{Ekh-JEMS} (cf. \cite{EES-TAMS},\cite{Ekholm-Honda-Kalman}, \cite{DRizell-lifting}) for more details and \cite{Chekanov}, \cite{Eliashberg-ICM98}, \cite{EGH}
for the original ideas underlying the construction.

Denote by $\cA (\Lambda,\lambda)$ a free non-commutative unital algebra over $\Z_2$ generated by
$\cR (\Lambda, \lambda)$.

The algebra $\cA (\Lambda,\lambda)$ comes with a filtration defined by the action: the action of a Reeb chord is its time-length and the action of a monomial which is a product of Reeb chords is the sum of the actions of the factors.
The action of the constant monomial $1\in\Z_2$ is set to be zero and the action of $0$ is defined as $-\infty$. For $r\in (-\infty,+\infty]$, define $\cA_r (\Lambda,\lambda)$ as the vector subspace of $\cA (\Lambda,\lambda)$ spanned over $\Z_2$ by the monomials whose action is smaller than $r$. It is easy to see that for a finite $r$ the vector space $\cA_r (\Lambda,\lambda)$ is finite-dimensional. For any $r\leq r'$,
there is a natural morphism $\cA_r (\Lambda,\lambda)\to \cA_{r'} (\Lambda,\lambda)$ induced by the inclusion of the spanning monomials.
For $r\leq 0$, $\cA_r (\Lambda,\lambda) := 0$.

For a suitable {\it cylindrical} (i.e. $\R_+$-equivariant)
almost complex structure $J$ on $\Sigma\times\R_+$ fitted to $s\lambda$,
one can define a
differential
$\partial_J$
on a generator of $\cA (\Lambda,\lambda)$
using a mod-$2$ count, up to domain reparameterization and the $\R_+$-action,
of rigid $J$-holomorphic maps of a disk with one positive and several (possibly no) negative boundary punctures into $\Sigma\times\R_+$. Such a map should send the boundary of the disk to $\Lambda\times\R_+$ and converge near each positive/negative puncture to a positive/negative cylinder over a Reeb chord in $\cR (\Lambda,\lambda)$. One also sets $\partial_J (1) :=0$. The differential is then extended by linearity and the Leibniz rule to the whole $\cA (\Lambda,\lambda)$.

Such a cylindrical $J$ for which all the moduli spaces of expected dimension $0$ (after
taking the quotient by the $\R_+$-action)
and relevant for the count above
are transversely cut-out and compact (meaning that $\partial_J$ is well-defined) and for which
also $\partial_J^2=0$, is called {\it regular} (with respect to $\Lambda\times \R_+$ and $\lambda$).
The set $\cJ (\Lambda,\lambda)$ of regular cylindrical $J$
is dense in the space of all cylindrical almost complex structures on $\Sigma\times \R_+$
\cite{Ekh-JEMS}, see also
\cite{EES-TAMS}, \cite{Ekholm-Honda-Kalman}, \cite{DRizell-lifting}.

Let $J\in \cJ (\Lambda,\lambda)$.

A computation using Stokes' theorem and similar to
\cite[Lemma 5.16]{BEHWZ} and \cite[Lemma B.3]{Ekh-JEMS} yields
\begin{equation}
\label{eqn-partial-J-respects-filtration}
\partial_J \left( \cA_r (\Lambda,\lambda)\right)\subset \cA_r (\Lambda,\lambda)\ \textrm{for all}\ r\in\R.
\end{equation}

Assume now that $\Lambda=\Lambda_0\sqcup\Lambda_1$ is a two-part Legendrian submanifold.


A monomial $a_1\cdot\ldots\cdot a_k$ in $\cA (\Lambda,\lambda)$ will be called {\it $ij$-composable}, or an
{\it $ij$-monomial}, for $i,j=0,1$,
if $a_k$ starts at $\Lambda_i$, $a_1$ ends at $\Lambda_j$, and for each $m=1,\ldots,k-1$ the beginning of $a_m$ lies in the same part of $\Lambda$ as the end of $a_{m+1}$. Note that, for $i\neq j$,
an $ij$-composable sequence of Reeb chords must contain at least one chord from $\cR_{ij} (\Lambda, \lambda)$.

Denote by $\cA (\Lambda_0,\Lambda_1,\lambda)$ the vector subspace of $\cA (\Lambda,\lambda)$ (over $\Z_2$) generated by all the $01$-composable monomials $a_1\cdot\ldots\cdot a_k$, $k\in\Z_{>0}$. (Note that the polynomials appearing in $\cA (\Lambda_0,\Lambda_1,\lambda)$ have no constant terms!). We will call $\cA (\Lambda_0,\Lambda_1,\lambda)$ the {\it $01$-subspace of $\cA (\Lambda,\lambda)$}.

For each $r\in (-\infty,+\infty]$,
set
\[
\cA_r (\Lambda_0,\Lambda_1,\lambda) := \cA (\Lambda_0,\Lambda_1,\lambda) \cap \cA_r (\Lambda,\lambda).
\]
In particular, $\cA_r (\Lambda_0,\Lambda_1,\lambda) = 0$ for $r\leq 0$.

Let $J\in \cJ (\Lambda=\Lambda_0\sqcup\Lambda_1,\lambda)$.

It is easy to see that $\cA (\Lambda_0,\Lambda_1,\lambda)$ is invariant under $\partial_J$ and therefore
\begin{equation}
\label{eqn-partial-J-respects-filtration-01}
\partial_J \left( \cA_r (\Lambda_0,\Lambda_1,\lambda)\right)\subset
\cA_r (\Lambda_0,\Lambda_1,\lambda)\ \textrm{for all}\ r\in\R.
\end{equation}

For each $r\in (-\infty,+\infty]$, define a vector space $V_r (\Lambda_0,\Lambda_1,\lambda,J)$ over $\Z_2$:
\[
V_r (\Lambda_0,\Lambda_1,\lambda,J) := \frac{\Ker \partial_J|_{\cA_r (\Lambda_0,\Lambda_1,\lambda)}}{\Image \partial_J|_{\cA_r (\Lambda_0,\Lambda_1,\lambda)}}.
\]
The inclusion maps $\cA_r (\Lambda_0,\Lambda_1,\lambda)\to \cA_{r'} (\Lambda_0,\Lambda_1,\lambda)$, $r\leq r'$, induce morphisms $V_r (\Lambda_0,\Lambda_1,\lambda,J)\to V_{r'} (\Lambda_0,\Lambda_1,\lambda,J)$.
For $r\leq 0$, we have $V_r (\Lambda_0,\Lambda_1,\lambda,J) = 0$.

It is easy to see that the vector spaces $V_r (\Lambda_0,\Lambda_1,\lambda,J)$, $r\in {(-\infty,+\infty]}$, together with the morphisms between them, form a persistence module supported in $(0,+\infty)$.
We will denote this persistence module by $V (\Lambda_0,\Lambda_1,\lambda,J)$.

Let $c>0$. Consider the contact form $c\lambda$ on $(\Sigma,\xi)$ and a regular cylindrical $J\in\cJ (\Lambda,\lambda)$.
(Note that $J$ is regular and cylindrical regardless of a choice of $c$).
Define $V (\Lambda_0,\Lambda_1,c\lambda,J)$
similarly to the above.
It is easy to see that
\[
V (\Lambda_0,\Lambda_1,c\lambda,J) = V^{[\times 1/c]} (\Lambda_0,\Lambda_1,\lambda,J).
\]


\hfill


\subsection{Morphisms of persistence modules defined by two-part-exact Lag\-ran\-gian co\-bor\-disms}
\label{subsec-morphisms-pers-mod-defd-by-Lagr-cobs}

With $(\Sigma,\xi)$ as above, assume, as in Section~\ref{sec-lch-pers-modules}, that

\begin{itemize}

\item{} $\lambda^-,\lambda^+$, $\lambda^-\prec\lambda^+$, are comfortable contact forms on $(\Sigma,\xi)$,

\item{} $(M,\theta)$ is a comfortable intergraphical exact symplectic cobordism, based on
$(\Sigma,\xi)$, between $(\Sigma^\pm,\lambda^\pm)$,

\item{} $L=L_0\sqcup L_1\subset (M,\theta)$ is a two-part-exact Lagrangian
cobordism between $(\Lambda^\pm,\lambda^\pm)$, where $\Lambda^\pm\subset (\Sigma,\xi)$ are
two-part Legendrian submanifolds,

\item{} the pairs $(\Lambda^\pm,\lambda^\pm)$ are non-degenerate (and thus, with our terminology,
$L$ is non-degenerate).

\end{itemize}

Let $S^+\cong\Sigma^+\times [1,+\infty)$
be a codimension-0 submanifold with boundary of $S\Sigma$ whose only boundary is
$\Sigma^+$ and that lies above $\Sigma^+$. Similarly,
let $S^-\cong\Sigma^-\times (0,1]$
be a codimension-0 submanifold with boundary of $S\Sigma$ whose only boundary is
$\Sigma^-$ and that lies below $\Sigma^-$.
Then we have
\[
S\Sigma = S^-\cup M\cup S^+.
\]
The 1-forms $s\lambda^\pm$ on $S^\pm$ can be glued with the 1-form $\theta$
on $M$ into a 1-form on $S\Sigma$ that we will denote by $\btheta$.
The 2-form $d\btheta$ is, of course, symplectic.

Set
\[
S\Lambda_{\lambda^\pm} := S\Lambda^\pm\cap S^\pm.
\]
These are Lagrangian submanifolds of $(S^\pm,d(s\lambda^\pm))$.
Define
\[
\bL :=  S\Lambda_{\lambda^-}\cup L \cup S\Lambda_{\lambda^+}.
\]
It is a Lagrangian submanifold of $(S\Sigma, d\btheta)$, called
{\it the completion of $L$}.

For an appropriate -- a so-called {\it adap\-ted (to $L$ and $\btheta$}) --
almost complex structure $I$ on $S\Sigma$
define a unital algebra morphism
\[
\Phi^{L,\theta,I}: \cA (\Lambda^+,\lambda^+) \to \cA (\Lambda^-,\lambda^-)
\]
by setting
\[
\Phi^{L,\theta,I} (1):=1
\]
and prescribing its values on the Reeb-chord generators:
for any $a\in \cR (\Lambda^+,\lambda^+)$,
\begin{equation}
\label{eqn-def-Phi-L-theta-I}
\Phi^{L,\theta,I} (a) :=
\sum\limits_{\dim \cM_{L,I} (a;b_1,\ldots,b_m) = 0} |\cM_{L,I} (a; b_1,\ldots,b_m)| b_1\cdot\ldots\cdot b_m,
\end{equation}
where $\cM_{L,I} (a;b_1,\ldots,b_m)$ is the moduli space, up to a domain reparameterization,
of $I$-holomorphic maps
of a disk with one positive and $m\geq 0$ negative boundary punctures into
$S\Sigma$ that send the boundary of the disk to $\bL$ and converge near
the positive puncture to a cylinder over $a$ and the negative punctures to the
cylinders over the chords $b_1,\ldots,b_m\in\cR (\Lambda^-,\lambda^-)$. (Note there might be repetitions
among $b_1,\ldots,b_m$, and in case $m=0$ there are no $b_i$'s).
The cardinality $|\cdot|$ of the moduli spaces is taken modulo $2$.
Here
$I^\pm$ are cylindrical almost complex structures on $S\Sigma$ induced
by the restrictions of $I$ to the ends of $S\Sigma$ -- we say that
$I^\pm$ are {\it cylindrical almost complex structures induced by $I$ on the ends}.
We say that such an $I$ is {\it regular} (with respect to $L$ and $\theta$)
if
\begin{itemize}

\item{} both $I^\pm$ are regular,

\item{} all the
moduli spaces in \eqref{eqn-def-Phi-L-theta-I} of expected dimension $0$ are compact and transversely cut out
(meaning that $\Phi^{L,\theta,I}$ is a well-defined unital
algebra morphism),

\item{} $\Phi^{L,\theta,I}\circ \partial_{I^+} = \partial_{I^-}\circ \Phi^{L,\theta,I}$.

\end{itemize}
See \cite{Ekh-JEMS} and \cite{Ekholm-Honda-Kalman} for more details.

The set $\cI (L,\theta)$ of regular $I$ is
dense in the space of all adapted almost complex structures
on $S\Sigma$ and the map
\[
\cI (L,\theta)\to \cJ (\Lambda^+,\lambda^+)\times \cJ (\Lambda^-,\lambda^-),\ I\mapsto (I^+,I^-),
\]
is surjective \cite{Ekh-JEMS}, cf. \cite{Ekholm-Honda-Kalman}.

Note that, with our notation, the results in \cite{Ekh-JEMS}, \cite{Ekholm-Honda-Kalman} are proved
in the setting where $\theta=\theta_{can}$. Clearly, the same proof works for any $\theta$ such that
$d\theta$ is identified
with $d(c\theta_{can})$, $c>0$, by a diffeomorphism of $M$ fixed near the boundary; this
assumption holds
by the definition of an intergraphical exact symplectic cobordism
(Definition~\ref{def-integraphical-exact-symplec-cob}).

With this terminology, a regular adapted almost complex structure $I$ on $(S\Sigma,\btheta)$ induces
regular cylindrical almost complex structures $I^\pm$ on $(S\Sigma,\btheta)$.

The map $\Phi^{L,\theta,I}$ is called {\it the cobordism map associated to $L,\theta,I$}.

Let us denote by $|\cdot|$ the actions of elements in the Chekanov-Eliashberg algebras.

For an $ij$-monomial
$a_1\cdot\ldots\cdot a_k$ in the Chekanov-Eliashberg algebra of a two-part Legendrian submanifold,
write
\[
\sigma (a_1\cdot\ldots\cdot a_k):= \sum_{l=1}^k \sigma (a_l) = j-i.
\]
For an empty monomial (i.e. $1$), we set $\sigma (1) :=0$.

The following result is based on an adaptation of
the Stokes computation in
 \cite[Lemma B.3]{Ekh-JEMS}.


\hfill


\begin{prop}[\cite{EP-big}]
\label{prop-action-inequality-cob-map}
If $\cM_{L,I} (a;b_1,\ldots,b_m)\neq\emptyset$ for some $a\in \cR (\Lambda^+,\lambda^+)$ and
$b_1,\ldots,b_m\in\cR (\Lambda^-,\lambda^-)$,
then
\[
\sum\limits_{i=1}^m |b_i| < |a| - \sigma(a)\gap_{L,\theta}.
\]
\end{prop}


\hfill


Since $I\in \cI (L,\theta)$, the cobordism map $\Phi^{L,\theta,I}$ is an algebra homomorphism,
and therefore
Proposition~\ref{prop-action-inequality-cob-map} yields the following corollary.


\hfill


\begin{cor}
\label{cor-action-inequality-cob-map}
Assume that $I\in \cI (L,\theta)$ and $a_1\cdot\ldots\cdot a_k\in \cA_r (\Lambda^+,\lambda^+)$ is
an $ij$-monomial, for $i,j=0,1$, $r\in\R$.

Then
\[
\Phi^{L,\theta,I} (a_1\cdot\ldots\cdot a_k)\in \cA_{r-(j-i)\gap_{L,\theta}} (\Lambda^-,\lambda^-).
\]
\Qed
\end{cor}


\hfill


Part A of the following proposition follows directly
from Corollary~\ref{cor-action-inequality-cob-map}.
In the case of a trivial comfortable exact symplectic cobordism,
this statement appeared in
\cite[Prop. 4.2]{EP-JFPTA2022}.
The proof of part B is straightforward.


\hfill


\begin{prop}
\label{prop-Phi-L-I-defines-morphism-of-pers-mod}
Let $I\in \cI (L,\theta)$ and
let $C:=G_{L,\theta}$.

Then:

\bigskip
\noindent
A. The cobordism map
$\Phi^{L,\theta,I}$ maps $\cA_r (\Lambda^+_0,\Lambda^+_1,\lambda^+)$
into $\cA_{r-C} (\Lambda^-_0,\Lambda^-_1,\lambda^-)$
for each $r\in (-\infty,+\infty]$
and, accordingly, defines a morphism of additive degree $C$ of persistence modules:
\[
\Phi^{L,\theta,I}_* : V (\Lambda^+_0,\Lambda^+_1,\lambda^+,I^+) \to V (\Lambda^-_0,\Lambda^-_1,\lambda^-,I^-).
\]

\bigskip
\noindent
B. For any $c>0$, we have $I\in \cI (L,c\theta)$ and
\[
\Phi^{L,c\theta,I}_* = \left(\Phi^{L,\theta,I}_*\right)^{[\times 1/c]}.
\]

\end{prop}


\hfill


The following proposition was stated in \cite[Prop. 4.4]{EP-JFPTA2022} and is proved in \cite{EP-big} using the well-known results
(see \cite{Ekh-JEMS}) about the cobordism map associated to a trivial exact Lagrangian cobordism.


\hfill


\begin{prop}[\cite{EP-big}, cf. \cite{EP-JFPTA2022}, Prop. 4.4]
\label{prop-Phi-L-trivial-cobordisms}
Let $0<s_-<s_+$.
Assume that the non-degenerate two-part-exact Lagrangian cobordism $L$ is trivial:
$L=(\Lambda_0\sqcup \Lambda_1)\times [s_-,s_+]\subset (\Sigma\times [s_-,s_+],\theta_{can})$.
Consider a regular cylindrical almost complex structure $J$ on $(\Sigma\times\R_+,s\lambda)$.
Set
\[
U:= V (\Lambda_0,\Lambda_1,s_+\lambda, J).
\]

Then
\[
\Phi^{L,s\lambda,J}_* = \Sh_U [\times s_+/s_-].
\]
\Qed
\end{prop}


\hfill


The following proposition is central for our considerations. It was stated in our paper
\cite[Prop. 4.5]{EP-JFPTA2022}
in the case of a trivial ambient exact symplectic cobordism and a constant family of cobordism forms.
Unfortunately,
the proof there was imprecise. We take this opportunity to state the proposition for
a larger class of exact symplectic cobordisms and two-part-exact Lagrangian isotopies
and present its corrected proof.


\hfill


\begin{prop}[cf. \cite{EP-JFPTA2022}, Prop. 4.5]
\label{prop-Phi-L-two-part-Lagr-isotopy}
With  $L\subset (M,\theta)$ being a non-degenerate
two-part-exact Lagrangian cobordism between $\Lambda^\pm = \Lambda^\pm_0\sqcup\Lambda^\pm_1$
as above, suppose that
$\{ L^\tau,\theta^\tau\}_{0\leq \tau\leq T}$, $L^0=L$, $\theta^0=\theta$, is a
two-part-exact Lagrangian cobordism isotopy fixed near the boundary. (In particular,
all $L^\tau$ are non-degenerate). Let $\gap (\tau)$ be the gap function of the isotopy.
Set:
\[
\gap_{min} := \min_{\tau\in [0,T]} \gap (\tau),
\]
\[
C_0 := \gap (0) - \gap_{min},\ C_T := \gap (T) - \gap_{min}.
\]

Let $I^0\in \cI (L^0,\theta^0)$, $I^T\in \cI (L^T,\theta^T)$, so that $I^0$, $I^T$ induce the same regular
cylindrical almost complex structures $J^\pm$ on $(\Sigma\times\R_+,s\lambda^\pm)$.

Then the following morphisms of persistence modules are equal:
\begin{equation}
\label{eqn-Phi-L0-LT-shifts}
\Sh (+C_0) \circ \Phi^{L^0,\theta^0,I^0}_* =
\Sh (+C_T) \circ \Phi^{L^T,\theta^T,I^T}_*.
\end{equation}
\end{prop}

\begin{proof}
For brevity, set
\[
V^\pm := V (\Lambda_0^\pm,\Lambda_1^\pm,\lambda^\pm,J^\pm).
\]
Denote the morphisms in the persistence modules $V^\pm$ respectively by $\pi^\pm_{s,t}$, $-\infty<s\leq t\leq +\infty$.

Let us fix $r\in\R$ such that
\[
r - \gap_{min}\notin \spec (V^-).
\]
By the discreteness of the spectrum in the definition of a persistence module
(Definition~\ref{def-persistence-module}), there exists $\varepsilon>0$ such that
\begin{equation}
\label{eqn-V-minus-semicontinuity}
[r-\gap_{min}, r-\gap_{min}+2\varepsilon]\cap \spec (V^-) = \emptyset.
\end{equation}

Let us fix this $\varepsilon>0$.

Consider a family of almost complex structures $\{ I^\tau\}_{0\leq \tau\leq T}$ on $S\Sigma$ connecting $I^0$ and $I^T$
so that each $I^\tau$ is adapted to $\theta^\tau$ and induces the same $J^\pm$ on the ends. It follows from \cite[Lemma B.8]{Ekh-JEMS} and \cite[Lemma 3.14]{Ekholm-Honda-Kalman} (which is a version of \cite[Lemma B.15]{Ekh-JEMS}) that if the family $\{ I^\tau\}$ is generic,
then there exist exceptional moments
\[
0:=\tau_0<\tau_1<\ldots<\tau_l<\tau_{l+1}:=T,
\]
such that the following claims hold:
\begin{itemize}

\item{} For any $\tau\in [0,T]\setminus \{ \tau_1,\ldots,\tau_l\}$, we have $I^\tau\in \cI (L^\tau,\theta^\tau)$
and therefore the cobordism map $\Phi^{L^\tau,\theta^\tau,I^\tau}$ is well-defined. We abbreviate
\[
\Phi^\tau := \Phi^{L^\tau,\theta^\tau,I^\tau}.
\]

\item{} The family $\{ \Phi^\tau\}$ of cobordism
maps undergoes changes only at $\tau_1,\ldots,\tau_l$.

Roughly speaking, as $\tau$ approaches an exceptional parameter $\tau_i$,
the index-$0$ pseudo-holomorphic disks in
$\cM_{L^\tau,I^\tau}(a;b_1,\ldots,b_m)$ may undergo a birth--death
or degenerate to broken configurations containing an index-$(-1)$
cobordism disk.
A birth--death creates or removes two disks with the same asymptotic
data, so the mod-$2$ cobordism counts remain unchanged;
at such a parameter we take the chain homotopy morphism $K_i$, introduced in formula \eqref{eqn-K-on-a-monomial-2} below, to be zero.

The appearance of index-$(-1)$ cobordism disks is analogous to what
happens in a homotopy of Morse continuation data.
Continuation trajectories of Fredholm index $-1$ may occur at isolated
parameter values, and counting these trajectories over the parameter
defines the Morse chain homotopy.

\item{} For each $i=1,\ldots,l$, and each sufficiently small $\delta>0$, such that $\tau_{i-1}+\delta\leq \tau_i - \delta
< \tau_i+\delta\leq \tau_{i+1} - \delta$,
the change from $\Phi^{\tau_i-\delta}$
to $\Phi^{\tau_i+\delta}$
can be described by the following formula:
\begin{equation}
\label{eqn-K-i-chain-homotopy}
K_i \circ \partial_{J^+} + \partial_{J^-}\circ K_i =
\Phi^{\tau_i+\delta}
- \Phi^{\tau_i-\delta},
\end{equation}
for a $\Z_2$-linear map
\[
K_i: \cA (\Lambda^+,\lambda^+)\to \cA (\Lambda^-,\lambda^-),\ K_i (1) = 0,
\]
that depends on $\delta>0$ and some additional data and is discussed below.

\item{} For each $i=1,\ldots,l$,
the map $K_i$ is defined on each monomial $a_1\cdot\ldots\cdot a_k$, $a_1,\ldots,a_k\in \cR (\Lambda^+,\lambda^+)$, as
\begin{equation}
\label{eqn-K-on-a-monomial-1}
K_i (a_1\cdot\ldots\cdot a_k) :=
\sum_{j=1}^k  \Phi^{\tau_i+\delta} (a_1\cdot\ldots\cdot a_{j-1})
\cdot K_i (a_j) \Phi^{\tau_i-\delta} (a_{j+1}\cdot\ldots\cdot a_k),
\end{equation}
so that for each $a\in \cR (\Lambda^+,\lambda^+)$
\begin{equation}
\label{eqn-K-on-a-monomial-2}
K_i (a):= \sum_{m\geq 0, b_1,\ldots, b_m\in \cR (\Lambda^-,\lambda^-)} n_{\tau_i} (a;b_1,\ldots,b_m) b_1\cdot\ldots\cdot b_m.
\end{equation}
Roughly speaking,
the coefficient $n_{\tau_i}(a;b_1,\ldots,b_m)$  counts modulo~$2$
(using a perturbation scheme involving additional choices
beyond the almost complex structures $\{I^\tau\}_{\tau_i-\delta\leq\tau\leq\tau_i+\delta}$, see \cite{Ekh-JEMS} for details)
the isolated pairs $(\tau,u_\tau)$, where

\begin{itemize}

\item{}
$\tau$ lies near $\tau_i$,

\item{}
$u_\tau$ is an index-$(-1)$ ``broken disk", or in more precise terms,
a configuration of disks with boundary punctures,
where each disk in the configuration
is a solution of a perturbed Cauchy-Riemann equation and where
the positive asymptotic of the whole configuration is defined by $a$ and the ordered negative asymptotics
are defined by $b_1,\ldots,b_m$.

\end{itemize}



Applying a
Stokes-type action inequality as in Proposition~\ref{prop-action-inequality-cob-map} to each component
in the configuration forming $u_\tau$ and then summing up these inequalities yields the inequality
\[
\sum\limits_{i=1}^m |b_i| < |a| - \sigma(a)\gap (\tau).
\]
A standard compactness argument, applied as the perturbations and the parameter windows around $\tau_i$
shrink to zero, yields then
\begin{equation}
\label{eqn-action-estimate-broken-minus-1-disk}
\sum\limits_{i=1}^m |b_i| \leq |a| - \sigma(a)\gap (\tau_i).
\end{equation}



\end{itemize}

Note that, with our notation, the results in \cite{Ekh-JEMS}, \cite{Ekholm-Honda-Kalman} are proved
in the setting where $\theta^\tau\equiv\theta_{can}$. One easily verifies
that the same proof works also if
the forms $d\theta^\tau$, $\tau\in [0,T]$, are identified
with $d(c\theta_{can})$, $c>0$, by diffeomorphisms of $M$ fixed on a common
collar neighborhood of the
boundary and depending smoothly on $\tau$;
this condition holds by the definition of an intergraphical exact symplectic cobordism
(Definition~\ref{def-integraphical-exact-symplec-cob}) and by the definition of a two-part-exact Lagrangian cobordism
isotopy in Section~\ref{subsec-leg-submfds-lagr-cobs}.

Assume that $\delta>0$ above is so small that
\begin{equation}
\label{eqn-gap-tau-i-tau-i-pm-delta}
|\gap (\tau_i) - \gap (\tau_i\pm\delta)|< \varepsilon\ \textrm{for all}\ i=1,\ldots,l.
\end{equation}

By \eqref{eqn-K-i-chain-homotopy}, we have
\begin{equation}
\label{eqn-sum-of-chain-homotopies}
\Phi^T
- \Phi^0
=
\sum_{i=1}^l \left(\Phi^{\tau_i+\delta}
- \Phi^{\tau_i-\delta}\right) =\\
\left(\sum_{i=1}^l K_i\right) \circ \partial_{J^+} + \partial_{J^-}\circ \left(\sum_{i=1}^l K_i\right).
\end{equation}

Let $i=1,\ldots,l$.

By \eqref{eqn-K-on-a-monomial-2},
\eqref{eqn-action-estimate-broken-minus-1-disk}
and Proposition~\ref{prop-action-inequality-cob-map},
for each $a\in \cR (\Lambda^+,\lambda^+)$,
\begin{equation}
\label{eqn-K-i-a}
|K_i (a)|\leq |a|-\sigma(a)\gap (\tau_i).
\end{equation}

Assume that $a_1\cdot\ldots\cdot a_k\in \cA (\Lambda^+_0,\Lambda^+_1,\lambda^+)$ is a $01$-monomial.
Consider any of the terms
\[
\Phi^{\tau_i+\delta} (a_1\cdot\ldots\cdot a_{j-1}) K_i (a_j)
\Phi^{\tau_i-\delta} (a_{j+1}\cdot\ldots\cdot a_k),\ j=1,\ldots,k,
\]
in the right-hand side of \eqref{eqn-K-on-a-monomial-1}. If $a_j\in \cR (\Lambda^+,\lambda^+)$
is an $i'j'$-chord, $i',j'=0,1$, then $a_1\cdot\ldots\cdot a_{j-1}$ is a $j'1$-monomial and
$a_{j+1}\cdot\ldots\cdot a_k$ is a $0i'$-monomial.
Corollary~\ref{cor-action-inequality-cob-map} and \eqref{eqn-gap-tau-i-tau-i-pm-delta} yield
\begin{equation}
\label{eqn-Phi-tau-i-plus-epsilon-action-estimate}
\begin{gathered}
\left|\Phi^{\tau_i+\delta} (a_1\cdot\ldots\cdot a_{j-1})\right|
\leq\\
\left| a_1\cdot\ldots\cdot a_{j-1}\right| - (1-j') \gap (\tau_i+\delta)<\\
\left| a_1\cdot\ldots\cdot a_{j-1}\right| - (1-j') \gap (\tau_i) +\varepsilon,
\end{gathered}
\end{equation}
\begin{equation}
\label{eqn-Phi-tau-i-minus-epsilon-action-estimate}
\begin{gathered}
\left|\Phi^{\tau_i-\delta} (a_{j+1}\cdot\ldots\cdot a_k)\right|
\leq\\
\left| a_{j+1}\cdot\ldots\cdot a_k\right| - i' \gap (\tau_i-\delta)<\\
\left| a_{j+1}\cdot\ldots\cdot a_k\right| - i' \gap (\tau_i) + \varepsilon.
\end{gathered}
\end{equation}
Combining \eqref{eqn-K-i-a}, \eqref{eqn-Phi-tau-i-plus-epsilon-action-estimate} and
\eqref{eqn-Phi-tau-i-minus-epsilon-action-estimate},  we get
\[
\left| \Phi^{\tau_i+\delta} (a_1\cdot\ldots\cdot a_{j-1}) K_i (a_j)
\Phi^{\tau_i-\delta} (a_{j+1}\cdot\ldots\cdot a_k)\right|<
\]
\[
<
\left| a_1\cdot\ldots\cdot a_{j-1}\right| -  (1-j')\gap (\tau_i) + \varepsilon +
\]
\[
+
|a_j| - (j'-i')\gap (\tau_i) +
\]
\[
+ \left| a_{j+1}\cdot\ldots\cdot a_k\right| - i' \gap (\tau_i) + \varepsilon =
\]
\[
=
\left| a_1\cdot\ldots\cdot a_k\right| - \gap (\tau_i) + 2\varepsilon \leq
\]
\[
\leq \left| a_1\cdot\ldots\cdot a_k\right| - \gap_{min} + 2\varepsilon.
\]
Since this holds for any $01$-monomial
$a_1\cdot\ldots\cdot a_k\in \cA (\Lambda^+_0,\Lambda^+_1,\lambda^+)$,
the formula \eqref{eqn-K-on-a-monomial-1} implies that for each $r\in\R$,
\[
K_i\left(\cA_r (\Lambda_0^+,\Lambda_1^+,\lambda^+)\right)\subset\cA_{r-\gap_{min}+2\varepsilon} (\Lambda_0^-,\Lambda_1^-,\lambda^-).
\]
Since this holds for all $i=1,\ldots,l$, using it together with
\eqref{eqn-partial-J-respects-filtration-01}
and
\eqref{eqn-sum-of-chain-homotopies} yields that
\[
\left(\Phi^T
- \Phi^0\right)
\left(\cA_r (\Lambda_0^+,\Lambda_1^+,\lambda^+)\right)
\subset\cA_{r-\gap_{min}+2\varepsilon} (\Lambda_0^-,\Lambda_1^-,\lambda^-),
\]
and
$\Sh_{V^-} (+C_T+2\varepsilon)\circ\Phi^T_*$ and $\Sh_{V^-} (+C_0+2\varepsilon)\circ\Phi^0_*$ induce the same morphism
\[
V^+_r\to V^-_{r-\gap_{min}+2\varepsilon}.
\]
Recalling \eqref{eqn-V-minus-semicontinuity}, we get that, in fact,
$\Sh_{V^-} (+C_T)\circ\Phi^T_*$ and $\Sh_{V^-} (+C_0)\circ\Phi^0_*$  induce the same morphism
\[
V^+_r\to V^-_{r-\gap_{min}},
\]
or, in other words,
\[
\left(\Sh_{V^-} (+C_T)\circ\Phi^T_*\right)_r = \left(\Sh_{V^-} (+C_0)\circ\Phi^0_*\right)_r.
\]
This equality holds for any $r\in\R$ such that $r-\gap_{min}\notin\spec (V^-)$.
Together with the discreteness, the semicontinuity
and the direct limit properties
in the definition of a persistence module (Definition~\ref{def-persistence-module}) this easily yields
the same equality for all $r\in (-\infty,+\infty]$.

This proves \eqref{eqn-Phi-L0-LT-shifts}.
%
\end{proof}

\hfill


%

The following proposition is proved in \cite{EP-big} using
the standard gluing techniques for pseudo-holomorphic curves with strip-like ends
 -- see \cite[Lemma 3.14]{Ekh-JEMS}, \cite[Cor. 3.11, Lemma 3.13]{Ekholm-Honda-Kalman},
as well as \cite{Chanda-AdvMath2025}.


\hfill


\begin{prop}[\cite{EP-big}]
\label{prop-concatenation-Lagr-cobs-composition-of-morphisms}
Assume that
$(M,\theta)$, $(\tM,\ttheta)$
are comfortable intergraphical exact symplectic cobordisms based on $(\Sigma,\xi)$
and
$L\subset (M,\theta)$, $\tL\subset (\tM,\ttheta)$
are non-degenerate two-part-exact Lagrangian cobordisms
such that the concatenation $L\sharp\tL\subset (M\sharp \tM, \theta\sharp \ttheta)$
is well-defined.

Let $I\in \cI (L,\theta)$, $\tI\in \cI (\tL,\ttheta)$, so that $I$ and $\tI$ induce the same regular
cylindrical almost complex structure, respectively, on the positive end of the completion of $(M,\theta)$ and on the negative end of
the completion of $(\tM,\ttheta)$.

Then there exists a regular adapted almost complex structure $I'\in \cI (L\sharp \tL, \theta\sharp \ttheta)$
such that
\[
\Phi^{L\sharp \tL,\theta\sharp\ttheta,I'}_* = \Phi^{L,\theta,I}_* \circ \Phi^{\tL,\ttheta,\tI}_*.
\]
\Qed
\end{prop}


\hfill


\subsection{The dependence on almost complex structures}
\label{subsec-dependence-on-acs}

\begin{prop}[\cite{EP-big}]
\label{prop-pers-module-indep-of-acs-up-to-isomorphism}
For a non-degenerate $\Lambda=\Lambda_0\sqcup\Lambda_1$ as above, different choices
of $J\in \cJ (\Lambda,\lambda)$ yield isomorphic persistence modules $V (\Lambda_0,\Lambda_1,\lambda,J)$.
\Qed
\end{prop}


\hfill


This justifies the notation $V (\Lambda_0,\Lambda_1,\lambda)$
used in Section~\ref{sec-lch-pers-modules} for the isomorphism class of the persistence modules
$V (\Lambda_0,\Lambda_1,\lambda,J)$ and the notation $\cB (\Lambda_0,\Lambda_1,\lambda)$ for
the corresponding barcode.

We will also write $\dim V_\infty (\Lambda_0,\Lambda_1,\lambda)$ for $\dim V_\infty (\Lambda_0,\Lambda_1,\lambda,J)$ --
which, by Proposition~\ref{prop-pers-module-indep-of-acs-up-to-isomorphism}, is independent of $J$.

The following proposition is a version of \cite[Prop. 2.19]{Ng2003}, where the result is
stated only for Legendrian links in the standard contact $\R^3$; for a proof in full generality, see \cite{EP-big}.


\hfill


\begin{prop}[\cite{Ng2003}; cf. \cite{EP-big}]
\label{prop-dim-V-infty-indep-of-Leg-isotopy}
The dimension $\dim V_\infty (\Lambda_0,\Lambda_1,\lambda)$ does not change under Legendrian isotopies
of $\Lambda=\Lambda_0\sqcup\Lambda_1$.
\end{prop}


\hfill


\begin{prop}[\cite{EP-big}]
\label{prop-Phi-L-I-indep-of-I-up-to-equivalence}
For a non-degenerate two-part-exact Lagrangian co\-bor\-dism $L\subset (M,\theta)$ as above,
the morphism
$\Phi^{L,\theta,I}_*$ is independent of $I\in \cI (L,\theta)$ up to the RL-equivalence of morphisms
in the category of persistence modules. The isomorphisms providing the equivalence are the ones supplied by
Proposition~\ref{prop-pers-module-indep-of-acs-up-to-isomorphism}.
\Qed
\end{prop}


\hfill


This justifies the notation $\Phi^{L,\theta}_*$
used in Section~\ref{sec-lch-pers-modules} for the RL-equivalence class of the morphisms
$\Phi^{L,\theta,I}_*$.


\hfill


\subsection{The retreat of two-part-exact Lagrangian cobordisms}
\label{subsec-retreat-of-Lagr-cobs}

Let us prove the theorems from Sections~\ref{subsec-main-results-on-retreat} and
\ref{subsec-appl-restriction-on-gap-fcns}.


\hfill


\noindent
{\bf Proof of Theorem~\ref{thm-trivial-morphism}:}
Follows from Propositions~\ref{prop-Phi-L-trivial-cobordisms} and
\ref{prop-retreat-of-mult-shift-Sh}.
\Qed


\hfill


\noindent
{\bf Proof of Theorem~\ref{thm-Ups-concatenation-of-cobordisms}:}
Follows from Proposition~\ref{prop-concatenation-Lagr-cobs-composition-of-morphisms}
and part B of Proposition~\ref{prop-Upsilon-composition-of-morphisms}.
\Qed


\hfill


\noindent
{\bf Proof of Theorem~\ref{thm-Ups-of-cobordism-positive}:}
Let us prove part A of the theorem.

By the definition, any two-part-exact Lagrangian cobordism
$L\subset (M=S\Sigma_{\mu^-,\mu^+},\theta)$ can be represented, for some sufficiently small
$\delta>0$, as
\[
L= K\sharp L'\sharp N
\]
for trivial two-part-exact Lagrangian cobordisms
\[
K\subset (S\Sigma_{\mu^-, (1+\delta)\mu^-}=\Sigma^-\times [1,1+\delta],s\lambda^-),
\]
\[
N\subset (S\Sigma_{(1-\delta)\mu^+,\mu^+} = \Sigma^+\times [1-\delta,1],s\lambda^+),
\]
and a two-part-exact Lagrangian cobordism
\[
L'\subset (S\Sigma_{(1+\delta)\mu^-, (1-\delta)\mu^+},\theta).
\]
Now part A of the theorem follows from Theorems~\ref{thm-trivial-morphism} and \ref{thm-Ups-concatenation-of-cobordisms}.

Part B follows from part B of Proposition~\ref{prop-Phi-L-I-defines-morphism-of-pers-mod}
and Proposition~\ref{prop-retreat-mult-shift-of-a-morphism}.

This finishes the proof.
\Qed


\hfill


\noindent
{\bf Proof of Theorem~\ref{thm-two-part-exact-Lagr-isotopy-morphisms-of-pmodules}:}
Proposition~\ref{prop-Phi-L-two-part-Lagr-isotopy}, together with Remark~\ref{rem-add-shift-RL-equivalence},
yields the following equality between the RL-equivalence classes of morphisms of persistence modules:
\begin{equation}
\label{eqn-Sh-C-0-Phi-L-0-Sh-C-T-Phi-L-T}
\Sh (+C_0) \circ \Phi^{L^0,\theta^0}_* =
\Sh (+C_T) \circ \Phi^{L^T,\theta^T}_*.
\end{equation}
Consequently,
\begin{equation}
\label{eqn-Ups-Sh-C-0-Phi-L-0-Ups-Sh-C-T-Phi-L-T}
\Ups \left(\Sh (+C_0) \circ \Phi^{L^0,\theta^0}_* \right) =
\Ups \left(\Sh (+C_T) \circ \Phi^{L^T,\theta^T}_*\right).
\end{equation}
By \eqref{eqn-Ups-Phi-leq-Ups-infty-Phi} and \eqref{eqn-Ups-infty-add-shift-composed-with-Phi-equals-Ups-infty-Phi},
\[
\Ups_\infty \left(\Sh (+C_0) \circ \Phi^{L^0,\theta^0}_* \right)\geq \Ups \left(\Sh (+C_0) \circ \Phi^{L^0,\theta^0}_* \right).
\]
By Theorem~\ref{thm-Ups-of-cobordism-positive},
Proposition~\ref{prop-retreat-of-add-shift-Sh} and part B of Proposition~\ref{prop-Upsilon-composition-of-morphisms},
\[
\Ups \left(\Sh (+C_T) \circ \Phi^{L^T,\theta^T}_*\right)\geq
\Ups \left(\Sh (+C_T)\right) + \Ups \left(\Phi^{L^T,\theta^T}_*\right)>
C_T.
\]
Combining the last two formulae with \eqref{eqn-Ups-Sh-C-0-Phi-L-0-Ups-Sh-C-T-Phi-L-T}, we get
\[
\Ups_\infty \left(\Phi^{L^0,\theta^0}_*\right) + C_0\geq
\Ups \left(\Sh (+C_0)\circ \Phi^{L^0,\theta^0}_*\right) \geq
\Ups \left(\Phi^{L^T,\theta^T}_* \right) + C_T,
\]
which is exactly the first claim of Theorem~\ref{thm-two-part-exact-Lagr-isotopy-morphisms-of-pmodules}.
The second one follows from
\eqref{eqn-Sh-C-0-Phi-L-0-Sh-C-T-Phi-L-T}, Proposition~\ref{prop-Ups-depends-only-on-isom-class-of-morphism}
and \eqref{eqn-Ups-infty-add-shift-composed-with-Phi-equals-Ups-infty-Phi}.

This finishes the proof.
\Qed


\hfill


\noindent
{\bf Proof of Theorem~\ref{thm-obstruction-to-Lagr-isotopies}:}
Suppose that the two-part-exact Lagrangian isotopy, fixed near the boundary, from
$L$ to $\tL$ has the form $\{ L^\tau\}_{0\leq\tau\leq T}$, $L=L^0$, $\tL = L^T$.
Let
\[
C_{min} := \min_{[0,T]} \gap_{L^\tau,s\lambda},\
C_0 := \gap_{L,s\lambda} - C_{min},\
C_T := \gap_{\tL,s\lambda} - C_{min}.
\]
Then, by Theorem~\ref{thm-two-part-exact-Lagr-isotopy-morphisms-of-pmodules},
\begin{equation}
\label{eqn-thm-two-part-exact-Lagr-isotopy-morphisms-of-pmodules}
\Ups_\infty \left(\Phi^{L,s\lambda}_*\right) + C_0\geq
\Ups \left(\Sh (+C_0)\circ \Phi^{L,s\lambda}_*\right) \geq
\Ups \left(\Phi^{\tL,s\lambda}_* \right) + C_T.
\end{equation}
Since
\[
\tL\cap \big(\Sigma\times [s_0-\delta,s_1+\delta]\big) =
\tLambda\times [s_0-\delta,s_1+\delta],
\]
$\tL$ can be represented as
\[
\tL = K\sharp \big(\tLambda\times [s_0,s_1]\big)\sharp N,
\]
for Lagrangian cobordisms $K=K_0\sqcup K_1\subset (\Sigma\times [s_-,s_0],s\lambda)$,
$N=N_0\sqcup N_1\subset (\Sigma\times [s_1,s_+],s\lambda)$.

The Lagrangian cobordisms $K,N$ are two-part-exact.
Indeed, since $\tL_0$, $\tL_1$, $\tLambda_0$, $\tLambda_1$
are connected,
so are $K_i$, $N_i$, $i=0,1$,
because we assume that all Lagrangian cobordisms have non-empty positive and negative boundaries.
For each $i=0,1$, choose a smooth function $f_i: \tL_i\to\R$ with
\[
df_i=s\lambda|_{\tL_i},
\]
\[
f_i=0\ \text{near the negative end of}\ \tL_i,
\]
\[
f_i=A_{\tL_i,s\lambda}\ \text{near the positive end of}\ \tL_i.
\]
Because $s\lambda$ vanishes on $\tLambda\times [s_0,s_1]$ and each $\tLambda_i$, $i=0,1$,
is connected, $f_i$ has one value $b_i$ on that entire cylinder. Then
\[
f_i|_{K_i}\ \text{and}\  f_i|_{N_i}-b_i
\]
are the required normalized primitives on each $K_i$, respectively $N_i$,
equal to $0$
near its negative boundary and to {\sl a single constant} near its positive boundary, so that
\[
A_{K_i,s\lambda}=b_i,\
A_{N_i,s\lambda}=A_{\tL_i,s\lambda}-b_i.
\]


Now, by Theorems~\ref{thm-Ups-of-cobordism-positive}, \ref{thm-trivial-morphism} and
\ref{thm-Ups-concatenation-of-cobordisms},
we get
\[
\Ups \left(\Phi^{\tL,s\lambda}_*\right) \geq
\Ups \left(\Phi^{K,s\lambda}_*\right)
+ \Ups \left(\Phi^{\tLambda\times [s_0,s_1],s\lambda}_*\right)
+ \Ups \left(\Phi^{N,s\lambda}_*\right) >
\]
\[
>
(s_1 - s_0) \lmintc \big(\cB (\tLambda_0,\tLambda_1,\lambda)\big).
\]
Combining it with \eqref{eqn-thm-two-part-exact-Lagr-isotopy-morphisms-of-pmodules}, we get
\[
\Ups_\infty \left(\Phi^{L,s\lambda}_*\right) + C_0\geq
\Ups \left(\Sh (+C_0)\circ \Phi^{L,s\lambda}_*\right) >
(s_1 - s_0) \lmintc \big(\cB (\tLambda_0,\tLambda_1,\lambda)\big) + C_T.
\]
In view of the definitions of $C_0$ and $C_T$, this yields
\[
\gap_{\tL,s\lambda} <
\gap_{L,s\lambda} +
\Ups_\infty\left(\Phi^{L,s\lambda}_*\right) -
(s_1-s_0) \lmintc \big(\cB (\tLambda_0,\tLambda_1,\lambda)\big),
\]
and if the Lagrangian isotopy $\{ L^\tau\}$ is gap-raising (meaning, in particular, that
$C_0 = 0$), then
\[
\gap_{\tL,s\lambda} <
\gap_{L,s\lambda} +
\Ups\left(\Phi^{L,s\lambda}_*\right) -
(s_1-s_0) \lmintc \big(\cB (\tLambda_0,\tLambda_1,\lambda)\big).
\]
This finishes the proof.
\Qed


\hfill


Let us recall and prove Theorem~\ref{thm-intro-application}.

Namely, let $Q$ be a closed connected manifold, $\dim Q >0$. Let $\Sigma := J^1 Q$ and let $\lambda$ be the standard contact 1-form
defining the standard contact structure on $\Sigma$.
By Example~\ref{exam-1-jet-space}, our previous results are applicable to $\Sigma$ and $\lambda$.

Let $F: Q\to\R$ be a positive Morse function and let $\kappa>0$.
Define
\[
a:= \min_Q F >0,
\]
\[
\Lambda^{\prime}:= j^1 0,
\]
\[
\Lambda^{\prime\prime} := j^1 F,
\]
\[
\Lambda^{\prime\prime\prime} := j^1 (F+\kappa).
\]

Clearly, $\Lambda:= \Lambda^{\prime}\sqcup \Lambda^{\prime\prime}$ and $\tLambda:=\Lambda^{\prime}\sqcup\Lambda^{\prime\prime\prime}$
are two-part Legendrian submanifolds of
$(\Sigma,\xi)$.
One easily checks that the pairs $(\Lambda,\lambda)$ and $(\tLambda,\lambda)$ are non-degenerate.

Define the trivial two-part-exact Lagrangian cobordism
\[
L:= (\Lambda^{\prime}\sqcup\Lambda^{\prime\prime})\times [s_-,s_+]\subset (\Sigma\times [s_-,s_+],s\lambda).
\]
It is non-degenerate.


\hfill


\begin{thm}
\label{thm-intro-application-copy}
Let
\[
0< s_- < s_0 < s_1 < s_+.
\]
Assume there exists an isotopy, fixed near the boundary, in the class of two-part-exact
Lagrangian cobordisms in $(\Sigma\times [s_-,s_+],s\lambda)$, from $L$ to a two-part-exact
Lagrangian cobordism $\tL = \tL_0\sqcup\tL_1\subset(\Sigma\times [s_-,s_+],s\lambda)$
whose intersection
with $\Sigma\times [s_0-\delta,s_1+\delta]$, for some $\delta>0$,
is the trivial Lagrangian cobordism
$(\Lambda^{\prime}\sqcup\Lambda^{\prime\prime\prime})\times [s_0-\delta,s_1+\delta]\subset (\Sigma\times [s_-,s_+],s\lambda)$.

Then
\begin{equation}
\label{eqn-bound-thm-intro-application-copy}
\gap_{\tL,s\lambda} < (s_+ - s_-) a -
(s_1 - s_0) (a+\kappa).
\end{equation}

%

\end{thm}


\begin{proof}
Let us start by recalling the following well-known observation.

For a positive Morse function $f: Q\to\R$, the
two-part Legendrian submanifold $\Lambda:=j^1 0\sqcup j^1 f$ is non-degenerate.
The Reeb chords of $\Lambda$
are all $01$-chords and the $01$-subspace $\cA (j^1 0, j^1 f,\lambda)$ is the linear span over $\Z_2$
of the $01$-chords of $\Lambda$.
These Reeb chords are also identified with the critical points of $f$; the action of a Reeb chord
is the corresponding critical value of $f$. For an appropriate regular cylindrical almost complex structure $J$,
the restriction of the differential $\partial_J$ to $\cA (j^1 0, j^1 f,\lambda)$
is identified with the differential
in the Lagrangian Floer complex of the projections of $j^1 0$ and $j^1 f$
to $T^* Q$ that are {\sl embedded} Lagrangian submanifolds -- see e.g. \cite{EGH} (the text following Prop. 2.8.2)
and \cite{Chantraine-et-al-pos-Leg-isotopies}
(the comment preceding Cor. 1.10).
The latter Floer differential is identified with
the Morse differential by the original work of Floer \cite{Floer-JDG88}.
Consequently, the barcode $\cB (j^1 0,j^1 f,\lambda)$ is the Morse barcode of $f$, i.e.
the barcode of
the persistence module defined by the Morse complex of $f$ filtered by the values of $f$ itself.
The Morse barcode of $f$ contains an infinite bar starting at $\min_Q f$ and this bar has the
smallest left end among all the bars in the barcode.

Set
\[
c:=\frac{s_+}{s_-},
\
\widetilde c:=\frac{s_1}{s_0}.
\]
Let us apply the observation above to the positive Morse functions $F$ and $F+\kappa$.
It yields
\begin{equation}
\label{eqn-lmintc-a-plus-kappa}
\lmintc \big(\cB (\Lambda^{\prime},\Lambda^{\prime\prime\prime},\lambda)\big) = a+\kappa,
\end{equation}
\begin{equation}
\label{eqn-lmintinfty-a}
\lminc \big(\cB (\Lambda^{\prime},\Lambda^{\prime\prime},\lambda)\big) = \lmininfty \big(\cB (\Lambda^{\prime},\Lambda^{\prime\prime},\lambda)\big) = a.
\end{equation}
By Theorem~\ref{thm-trivial-morphism},
\begin{equation}
\label{eqn-Ups-infty-lminc}
\Ups_\infty (\Phi^{L,s\lambda}_*) = (s_+ - s_-) \lmininfty \big(\cB (\Lambda^{\prime},\Lambda^{\prime\prime},\lambda)\big).
\end{equation}


Corollary~\ref{cor-obstruction-to-Lagr-isotopies-trivial-cobordism} is applicable in our setting.
Note that the hypothesis in Corollary~\ref{cor-obstruction-to-Lagr-isotopies-trivial-cobordism} (inherited from Theorem~\ref{thm-obstruction-to-Lagr-isotopies}) about the connectedness of
$\tLambda_i$, $L_i$, $\tL_i$, $i=0,1$, is satisfied since $Q$ is assumed to be connected.

Now, since $\gap_{L,s\lambda} =0$,
the equations \eqref{eqn-lmintc-a-plus-kappa}, \eqref{eqn-lmintinfty-a},
\eqref{eqn-Ups-infty-lminc}, together with
Corollary~\ref{cor-obstruction-to-Lagr-isotopies-trivial-cobordism},
yield
\[
\gap_{\tL,s\lambda} < (s_+ - s_-) a -
(s_1 - s_0) (a+\kappa),
\]
which is precisely \eqref{eqn-bound-thm-intro-application-copy}.

This finishes the proof.
\end{proof}


\hfill


\begin{exam}
\label{exam-estimate-sharp}
{\rm
Let us present an example in which the bound
\begin{equation}
\label{eqn-bound-c-a-tc-kappa-geq-0-copy}
G_{\tL,s\lambda} < (s_+ - s_-) a -
(s_1 - s_0) (a+\kappa),
\end{equation}
appearing in Theorem~\ref{thm-intro-application-copy}, is sharp.

Consider the setting of Theorem~\ref{thm-intro-application-copy}.
Assume that
\[
s_0 - s_- = s_1 - s_0 = s_+ - s_1 =:r
\]
and
\[
G_{\tL,s\lambda} = 0.
\]
In this case
\eqref{eqn-bound-c-a-tc-kappa-geq-0-copy} can be rewritten as
\[
\frac{\kappa}{a} < 2.
\]

In order to show that this bound obstructing the existence of $\tL$
is sharp, let us show that, for $s_0, s_1$ as above and for $\kappa/a$
arbitrarily close to $2$ from below, there does exist the required $\tL$.

Write
\[
L_0 := \Lambda^{\prime}\times [s_-,s_+],
\]
\[
L_1 := \Lambda^{\prime\prime}\times [s_-,s_+],
\]
so that
\[
L = L_0\sqcup L_1.
\]

Consider a symplectic Hamiltonian flow
$\{\phi_\tau\}$ on $J^1 Q\times\R_+$,
equal to the identity for all $\tau$
near $\Sigma^\pm = \Sigma\times \{ s_\pm\}$ and
generated by a Hamiltonian
$h(s)$ whose graph, in the case $h(s_+) = 0$,
is illustrated in Figure~\ref{fig-graph-of-hamiltonian}.



\hfill


\begin{figure}[htbp]
\centering
\includegraphics[width=0.7\textwidth]{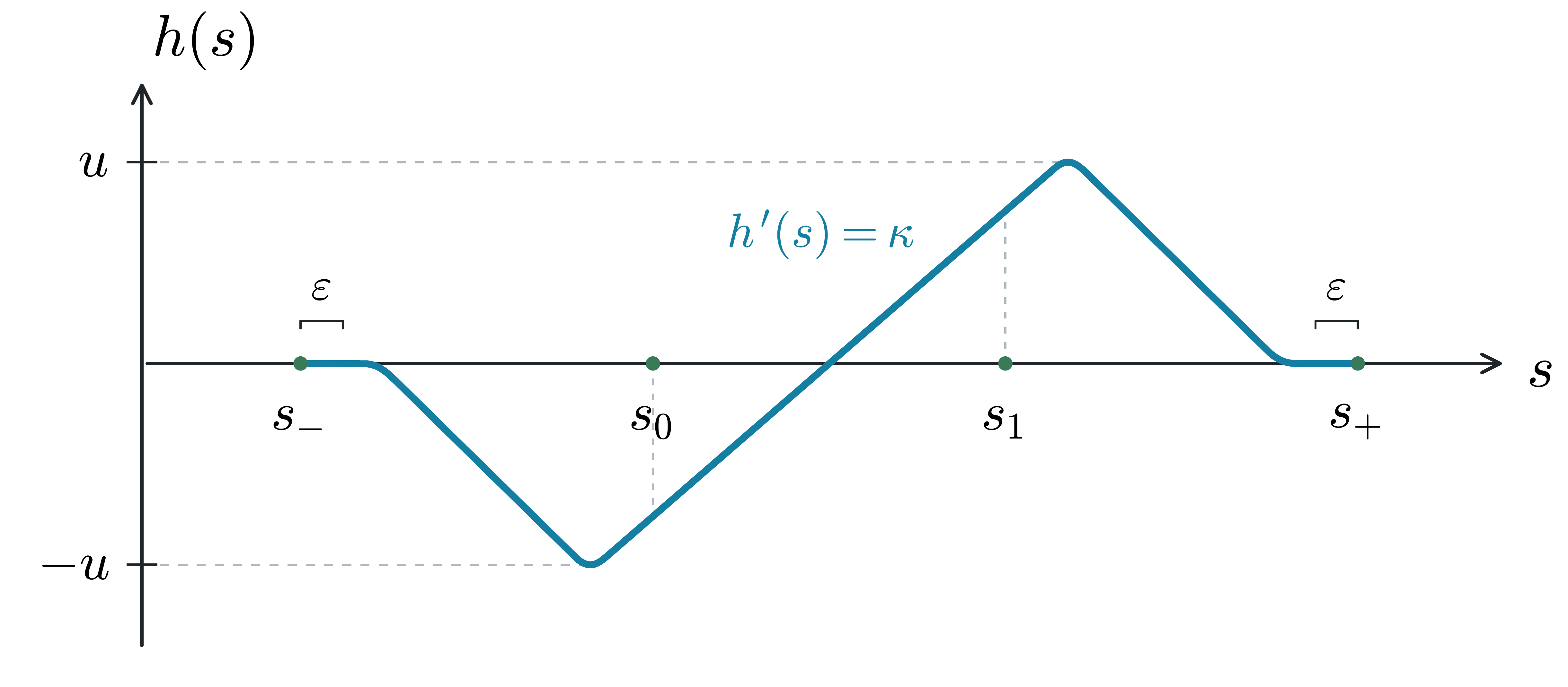}
\caption{}
\label{fig-graph-of-hamiltonian}
\end{figure}


\hfill


In particular, we assume the following:

\smallskip
\noindent
I. $h$ vanishes on an $\epsilon$-neighbourhood of $s_-$ and equals a constant $h(s_+)\geq 0$
on an $\epsilon$-neighbourhood of
$s_+$, for a small
$\epsilon>0$.

\smallskip
\noindent
II. $h$ is linearly increasing on an open interval containing $[s_0,s_1]$ with the slope $\kappa$.

\smallskip
\noindent
III.
$h'(s) > - a$ for all
\[
s \in [s_-,s_0] \cup [s_1,s_+].
\]

\smallskip
Let us discuss how to make sure that conditions II and III are both satisfied by $h$.

If we require $h(s_+)= 0$, then, for a sufficiently small $\epsilon >0$,
a function satisfying I-III can be constructed
whenever
\[
- \min_{\R_+} h = \max_{\R_+} h =: u > 0
\]
and
\begin{equation}
\label{eqn-u-geq-kappa-r-divided-by-two}
u>\kappa r/2,
\end{equation}
\begin{equation}
\label{eqn-u-smaller-a-r}
u < ar.
\end{equation}
The conditions \eqref{eqn-u-geq-kappa-r-divided-by-two} and \eqref{eqn-u-smaller-a-r}
are satisfied by the same $u$ if and only if
\[
\kappa/2 < a.
\]
In other words,
given $\kappa r/2<u<ar$, one can choose a smooth $h$ satisfying I-III with extrema $\pm u$.

If we do not require that $h(s_+) = 0$, then conditions II and III can be easily achieved for
$h$ with an arbitrarily large $h(s_+) > 0$.

Now let us write $J^1 Q \times \R_+ = T^* Q\times \R (z)\times \R_+ (s)$. Then we have
\begin{equation}
\label{eq-flow}
\phi_\tau (x,z,s) = (x, z+ h'(s)\tau, s),
\end{equation}
for any
$(x,z,s)\in T^* Q\times\R\times\R_+$.

By \eqref{eq-flow}, condition III
yields
\begin{equation}
\label{eq-int}
\phi_\tau (L_1) \cap L_0 =\emptyset \;\;\forall \tau \in [0,1]\;.
\end{equation}

Consequently, we
get a two-part-exact Lagrangian cobordism
isotopy, fixed near the boundary, of the form
$\{L^\tau = L^\tau_0\sqcup L^\tau_1\}_{0\leq\tau\leq 1}$, where
\[
L^\tau_0 = L_0,
\]
\[
L^\tau_1 = \phi_\tau (L_1)
\]
for all $\tau$
and
$L=L^0$.

Let
\[
\tL:= L^1.
\]

Condition II implies that, for a sufficiently small $\delta>0$ and $I_\delta := [s_0-\delta,s_1+\delta]$,
\begin{equation}
\label{eqn-tL-intersection-s-0-s-1}
\tL\cap (J^1 Q\times I_\delta) = \left(\Lambda^{\prime}\times I_\delta\right)\sqcup
\phi_1\left(\Lambda^{\prime\prime}\times I_\delta\right) =
\left(\Lambda^{\prime}\sqcup\Lambda^{\prime\prime\prime}\right)\times I_\delta.
\end{equation}

An easy computation shows that the gap function of the isotopy
$\{ L^\tau\}_{0\leq\tau\leq 1}$
is $-\tau h(s_+)$. Thus, if $h(s_+)=0$, the gap function is identically $0$ (which means, in particular,
that $G_{\tL,s\lambda} = 0$)
and this is a gap-raising isotopy.

Summing up, we get that:

\begin{itemize}

\item{}
As long as $h(s_+) = 0$ and $\kappa/a <2$ -- and, in particular,
if $\kappa/a$ is arbitrarily
close to $2$ from below -- we have a two-part-exact Lagrangian
cobordism $\tL\subset (J^1 Q\times [s_-,s_+], s\lambda)$
such that there exists a gap-raising Lagrangian isotopy from $L$
to $\tL$ and \eqref{eqn-tL-intersection-s-0-s-1} is satisfied. This means that
the constraint \eqref{eqn-bound-c-a-tc-kappa-geq-0-copy} is sharp in this example.

\item{} If we allow $h(s_+) >0$, then one can construct a two-part-exact Lagrangian
cobordism $\tL\subset (J^1 Q\times [s_-,s_+], s\lambda)$, so that:

\smallskip
\noindent
- there exists a two-part-exact Lagrangian isotopy from $L$ to $\tL$ fixed near the boundary,

\smallskip
\noindent
- $\tL$ satisfies \eqref{eqn-tL-intersection-s-0-s-1} and has
a negative gap arbitrarily close to $-\infty$.

\end{itemize}

}
\end{exam}


\hfill


\noindent
{\bf Acknowledgments:} We thank F.Bourgeois, B.Chantraine, T.Ekholm, G.Di\-mi\-tro\-glou Rizell and M.Sullivan
for helpful discussions over many years on mathematical questions related to this paper.




\hfill


\noindent
{\bf AI declaration:} In the preparation of this paper, we used the assistance of GPT-6 Astra in the following ways:

\begin{itemize}

\item{} It suggested the use of the $\Upsilon_\infty$-version of the retreat (as opposed to the usual retreat $\Ups$)
in the proof of Theorem~\ref{thm-intro-application}.

\item{} It suggested the second example in Remark~\ref{rem-examples-of-barcodes}.

\item{} It found a problem with the proof of \cite[Prop. 4.5]{EP-JFPTA2022} and suggested the correction.

\item{} Proofreading.

\end{itemize}


\bibliographystyle{alpha}

\noindent {\sc Michael Entov\\
Faculty of Mathematics,\\
Technion - Israel Institute of Technology,\\
Haifa 3200003, Israel}\\
{\tt  entov@technion.ac.il}
\\

\noindent {\sc Leonid Polterovich\\
School of Mathematical Sciences,\\
Tel Aviv University,\\
Tel Aviv 69978, Israel}\\
{\tt  polterov@tauex.tau.ac.il}
\\

\end{document}